\documentclass[11pt]{amsart}
\usepackage{amsxtra}
\usepackage{amssymb}
\usepackage[mathscr]{eucal}
\usepackage{musicography}
\theoremstyle{plain}

\newtheorem{theorem}{Theorem}[subsection]
\newtheorem{lemma}[theorem]{Lemma}
\newtheorem{definition-theorem}[theorem]{Definition-Theorem}
\newtheorem{proposition}[theorem]{Proposition}
\newtheorem{corollary}[theorem]{Corollary}
\newtheorem{definition}[theorem]{Definition}
\newtheorem{example}[theorem]{Example}
\newtheorem{remark}[theorem]{Remark}
\newtheorem{conjecture}[theorem]{Conjecture}
\newtheorem{notation}[theorem]{Notation}

\newtheorem{assumption}[theorem]{Assumption}
\newtheorem*{maintheorem*}{Main Theorem}
\newcommand \bth[1] { \begin{theorem}\label{t#1} }
\newcommand \ble[1] { \begin{lemma}\label{l#1} }

\newcommand \bpr[1] { \begin{proposition}\label{p#1} }
\newcommand \bco[1] { \begin{corollary}\label{c#1} }
\newcommand \bde[1] { \begin{definition}\label{d#1}\rm }
\newcommand \bex[1] { \begin{example}\label{e#1}\rm }
\newcommand \bre[1] { \begin{remark}\label{r#1}\rm }
\newcommand \bcj[1] { \begin{conjecture}\label{j#1}\rm }

\newcommand \bnota[1] { \begin{notation}\label{n#1}\rm }
\renewcommand {\eth} { \end{theorem} }
\newcommand {\ele} { \end{lemma} }

\newcommand {\epr} { \end{proposition} }
\newcommand {\eco} { \end{corollary} }
\newcommand {\ede} { \end{definition} }
\newcommand {\eex} { \end{example} }
\newcommand {\ere} { \end{remark} }
\newcommand {\ecj} { \end{conjecture} }

\newcommand {\enota} { \end{notation} }
\newcommand \thref[1]{Theorem \ref{t#1}}
\newcommand \leref[1]{Lemma \ref{l#1}}
\newcommand \prref[1]{Proposition \ref{p#1}}
\newcommand \coref[1]{Corollary \ref{c#1}}

\newcommand \deref[1]{Definition \ref{d#1}}

\makeatletter
\newsavebox{\@brx}
\newcommand{\llangle}[1][]{\savebox{\@brx}{\(\m@th{#1\langle}\)}%
  \mathopen{\copy\@brx\kern-0.5\wd\@brx\usebox{\@brx}}}
\newcommand{\rrangle}[1][]{\savebox{\@brx}{\(\m@th{#1\rangle}\)}%
  \mathclose{\copy\@brx\kern-0.5\wd\@brx\usebox{\@brx}}}
\makeatother

\usepackage{tikz}
\usetikzlibrary{matrix}
\usepackage{graphicx,ctable,booktabs}
\usetikzlibrary{arrows.meta}
\usepackage{tikz-cd}

\usepackage{stmaryrd}

\usepackage{hyperref}

 \DeclareMathOperator{\Proj}{Proj}

\DeclareMathOperator{\End}{End} 
 \DeclareMathOperator{\Hom}{Hom}

\DeclareMathOperator{\Mod}{{\sf Mod}}
\DeclareMathOperator{\modd}{{\sf mod}}

\DeclareMathOperator{\Spc}{Spc}

\DeclareMathOperator{\ev}{ev}
\DeclareMathOperator{\coev}{coev}

\DeclareMathOperator{\perf}{perf}
\DeclareMathOperator{\Ind}{Ind}

\newcommand{\ul}{\underline}

\newcommand{\mf}{\mathfrak}
\newcommand{\mc}{\mathcal}
\newcommand{\mb}{\mathbb}
\newcommand{\id}{\operatorname{id}}

\newcommand{\kk}{\Bbbk}
\newcommand{\idem}{\musNatural{}}
\newcommand{\bT}{\mathbf T}

\newcommand{\bS}{\mathbf S}
\newcommand{\bC}{\mathbf C}
\newcommand{\bK}{\mathbf K}
\newcommand{\bL}{\mathbf L}
\newcommand{\bP}{\mathbf P}
\newcommand{\bQ}{\mathbf Q}
\newcommand{\bI}{\mathbf I}
\newcommand{\bJ}{\mathbf J}
\newcommand{\bM}{\mathbf M}
\newcommand{\bN}{\mathbf N}

\newcommand{\Loc}{\operatorname{Loc}}

\newcommand{\unit}{\ensuremath{\mathbf 1}}

\newcounter{listequation}

\numberwithin{equation}{subsection}
\begin{document}
\title[Chinese Remainder Theorem and Carlson's Theorem]{A Chinese Remainder Theorem and \\ Carlson's Theorem for 
\\
Monoidal Triangulated Categories}
\author[Daniel K. Nakano]{Daniel K. Nakano}
\address{Department of Mathematics \\
University of Georgia \\
Athens, GA 30602\\
U.S.A.}
\email{nakano@math.uga.edu}
\thanks{The research of D.K.N was supported in part by NSF grant DMS-2101941.
The research of  K.B.V. was supported in part by NSF Postdoctoral Fellowship DMS-2103272.
The research of M.T.Y. was supported in part by NSF grant DMS–2200762.}
\author[Kent B. Vashaw]{Kent B. Vashaw}
\address{
Department of Mathematics \\
MIT \\
Cambridge, MA 02139 \\
U.S.A.
}
\email{kentv@mit.edu}
\author[Milen T. Yakimov]{Milen T. Yakimov}
\address{
Department of Mathematics \\
Northeastern University \\
Boston, MA 02115 \\
U.S.A.}
\email{m.yakimov@northeastern.edu}
\begin{abstract}
 In this paper the authors prove fundamental decomposition theorems pertaining to the internal structure of monoidal triangulated categories (M$\Delta$Cs). The tensor structure of an M$\Delta$C enables one to view these categories like (noncommutative) rings and to attempt to extend the key results for the latter to the categorical setting. The main theorem is an analogue of the Chinese Remainder Theorem involving the Verdier quotients for coprime thick ideals. 
This result is used to obtain orthogonal decompositions of the extended endomorphism rings of idempotent algebra objects of M$\Delta$Cs. The authors also provide topological characterizations on when an M$\Delta$C contains a pair of coprime proper thick ideals, and additionally, when the latter are complementary in the sense that their intersection is contained in the prime radical of the category.  

As an application of the aforementioned results, the authors establish for arbitrary M$\Delta$Cs a general version of Carlson's theorem on the connnectedness of supports 
for indecomposable objects. Examples of our results are given at the end of the paper for the derived category of schemes and for the stable module categories for finite 
group schemes. 
\end{abstract}
\maketitle
\section{Introduction}
\label{Intro}
\subsection{Monoidal triangulated categories}  Tensor triangular geometry (TTG) was introduced in the mid 2000's via tensor triangulated categories (TTC) by Balmer to provide a unifying method for understanding the underlying geometry 
that arises from examples in representation theory, homotopy theory, algebraic geometry, commutative algebra and category theory \cite{Balmer}. Given a TTC, $\bf K$, Balmer's idea was to utilize the (symmetric) tensor structure to view the category like a commutative ring, and construct a topological space $\Spc \bK$ (upgraded to a locally ringed space) now known at the Balmer spectrum via prime ideals in $\bf K$. 

The authors \cite{NVY1} developed a general noncommutative version of Balmer's theory that deals with an arbitrary monoidal triangulated category $\bK$ (M$\Delta$C). See also the work in \cite{BKS1}.  
There are many important families of M$\Delta$Cs that arise from diverse areas of mathematics and mathematical physics. One such family is comprised of the stable module categories of 
finite dimensional Hopf algebras that includes quantum groups, and more generally, the stable categories of finite tensor categories \cite{EGNO}.
These M$\Delta$Cs are in general not symmetric (because Hopf algebras are not necessarily cocommutative), and should be viewed like a noncommutative ring. 

The key feature of our new approach is to define the noncommutative Balmer spectrum $\Spc \bK$ and support data for $\bK$ via tensoring thick ideals of $\bK$, and not to use object-wise tensoring. 
In \cite{NVY2}, our approach enables us to distinguish between a prime spectrum, $\Spc \bf K$, and a completely prime spectrum, $\text{CPSpc} \bK$, and to provide intrinsic categorical answers to questions in the literature about the tensor product property on supports (cf. \cite[\S 3 and 4]{NVY2}), which were previously treated on a case-by-case basis of examples and counterexamples. 

\subsection{Results in the paper} In order to gain a better understanding about monoidal triangular geometry, the authors will treat M$\Delta$Cs like rings and develop the structure and representation theory. 
The first step in such a process entails developing an internal structure theory. For an arbitrary ring $R$, the Chinese Remainder Theorem states that if $I_1$ and $I_2$ are relatively prime ideals (i.e., $R=I_1+I_2$), and $I=I_1\cap I_2$ then $R/I=R/I_1 \oplus R/I_2$. See \cite[p 110, 10, 11]{Jac}. Our first main result is such a theorem for an arbitrary M$\Delta$C through the use of localization functors.
It can be thought of as analogous to block decompositions that encompass specific cases proved earlier in \cite{CDW1994}, \cite{BK2000}. Our results also provide a decomposition theorem for related indempotent algebra objects of an M$\Delta$C, which are algebra objects whose product map is an isomorphism, see Definition \ref{dalg}; these are the natural monoidal triangulated categorical analogues of idempotents of an algebra.   

Recall that an M$\Delta$C, $\bK$, is called small \cite{BS} it is isomorphic to the compact part of a compactly generated M$\Delta$C, $\widehat{\bK}$, called a big M$\Delta$C.
\medskip

\noindent{\bf{Theorem A.}} {\em{Let $\bK$ be a small M$\Delta$C which is the compact part of an associated big M$\Delta$C, $\widehat \bK$, and $\bI_1, \ldots, \bI_n$ be a collection of thick ideals of $\bK$ that are pairwise coprime, see \deref{def-ideals}(a).}} 
\begin{enumerate}
    \item[(a)] We have the equivalences of M$\Delta$Cs 
\begin{align*}
\widehat \bK/\Loc(\bI_1 \cap \ldots \cap \bI_n) &\cong \widehat \bK / \Loc( \bI_1) \oplus \ldots \oplus \widehat \bK / \Loc(\bI_n) \quad \mbox{and}
\\
\left ( \bK/(\bI_1 \cap \ldots \cap \bI_n) \right )^{\idem} &\cong (\bK / \bI_1)^{\idem} \oplus \ldots \oplus (\bK / \bI_n)^{\idem},
\end{align*}
where $\Loc$ refers to localizing ideal and $()^{\idem}$ to idempotent completion.
\item[(b)] {\em{The extended endomorphism ring \eqref{ext-end} and the Tate extended endomorphism ring \eqref{ext-endT} of the idempotent algebra object $L_{\bI_1 \cap \ldots \cap \bI_n}(\unit)$
of $\widehat{\bK}$:
\[
\End_{\widehat{\bK}}^\bullet
\big(L_{\bI_1 \cap \ldots \cap \bI_n}(\unit) \big) \cong \bigoplus_{j=1}^n \End_{\widehat{\bK}}^\bullet
\big(L_{\bI_j}(\unit) \big), 
\quad
\widehat{\End}_{\widehat{\bK}}^\bullet
\big( L_{\bI_1 \cap \ldots \cap \bI_n} (\unit) \big) \cong
\bigoplus_{j=1}^n 
\widehat{\End}_{\widehat{\bK}}^\bullet
\big(L_{\bI_j}(\unit) \big), 
\]
where $L_\bI$ refers to the localizing functor, see Sect. \ref{sec:loc-funct}.}}
\end{enumerate}
\medskip

The second part of the theorem is a categorical version of orthogonal idempotent decompositions in algebras; it decomposes the two forms of the extended endomorphism rings of the idempotent algebra object 
$L_{\bI_1 \cap \ldots \cap \bI_n}(\unit)$ in terms of the orthogonal
idempotent algebra objects $L_{\bI_1}(\unit), \ldots, L_{\bI_n}(\unit)$.

For commutative rings, the nilradical is an important ideal that can be characterized by the intersection of prime ideals or by taking the set of nilpotent elements in the ring (cf. \cite[Corollary 4.5]{Kunz}). There are various radical theories for noncommutative rings. The analogue of the nilradical is the Baer radical (also known as the prime radical or the lower nilradical). See the work of Amitsur \cite{Amit1, Amit2, Amit3}. The prime/Baer radical is defined by taking the intersection of all prime ideals; it is consists of nilpotent elements, and is a nilpotent ideal under a weak Noetherianity assumption \cite[Theorem 3.11]{GW}. In this paper we will use an analogue of the prime/Baer radical for M$\Delta$Cs to measure disjointness of thick ideals of M$\Delta$Cs. 

Our second theorem establishes topological characterizations of the existence of pairs of coprime proper thick ideals of M$\Delta$Cs.
\medskip

\noindent
{\bf{Theorem B.}}
\begin{enumerate}
\item[(a)] {\em{For every M$\Delta$C, $\bK$, the following are equivalent:}}
\begin{itemize}
    \item[(i)] {\em{There exist two coprime proper thick ideals of $\bK$.}}
    \item[(ii)] {\em{$\Spc \bK$ is reducible, that is, there exist proper closed sets $Z_1$ and $Z_2$ in $\Spc \bK$ with $\Spc \bK=Z_1 \cup Z_2.$}}
\end{itemize}
\item[(b)] {\em{If an M$\Delta$C, $\bK$, is generated by a single object as a thick subcategory or its Balmer spectrum $\Spc \bK$ is Noetherian then following are equivalent:}}
\begin{enumerate}
    \item[(i)] {\em{There exist two proper coprime thick ideals of $\bK$ that are complementary in the sense that their intersection is contained in the prime radical of $\bK$.}}
    \item[(ii)] {\em{$\Spc \bK$ is disconnected, that is, there exist proper disjoint closed sets $Z_1$ and $Z_2$ in $\Spc \bK$ with $\Spc \bK=Z_1 \sqcup Z_2.$}}
\end{enumerate}
\end{enumerate}
\medskip

The motivation for the term complementary in part (b) of the theorem is coming from the fact that, if each object of $\bK$ is left or right dualizable, then its prime radical is $\{0\}$, i.e., in this important and general situation the complementary condition for thick ideals means that their intersection is $\{0\}$. 

The special case of part (b) of Theorem B when $\bK$ is the stable module category of a finite group was proved by Krause in \cite[Theorem C]{Krause0}, stated there in an equivalent formulation using the cohomological support. 

Theorem A is used to prove a Carlson Connectedness Theorem \cite{Carlson1984} for arbitrary monoidal triangulated categories.
\medskip

\noindent
{\bf{Theorem C.}} {\em{Assume that $\bK$ is a small M$\Delta$C that is generated by a single object as a thick subcategory or that its Balmer spectrum $\Spc \bK$ is Noetherian.
Then for any indecomposable object $A \in \bK$, its support set $V(A)$ in the Balmer spectrum $\Spc \bK$ is connected.}}
\medskip

We view the results in this paper as foundational for future work in continued development of monoidal triangular geometry.  

\subsection{Organization of the paper.} In the next section (Section 2), we provide background on noncommutative tensor triangular geometry and localization functors in triangulated categories. The {\em prime radical} 
is also introduced. In Section 3, the Chinese Remainder Theorem \ref{tCRT} is stated and proved. Applications of this theorem are provided in Section 4 that yield orthogonal decompositions of the extended endomorphism rings of idempotent algebra objects (see \thref{end-algs}). 

In Section 5, we characterize M$\Delta$Cs that have at least two coprime proper ideals. 
This characterization is carried out further in Section 6 where the assumption that the ideals are complementary is added. Section 7 culminates in our general version 
of Carlson's Connectedness Theorem for M$\Delta$Cs. Examples of our results are presented in Section 8. Finally, we have included an Appendix that contains several auxiliary results that are of independent interests and which are used in the paper. 

\subsection*{Acknowledgements}
The authors are grateful to Pavel Etingof for extensive discussions. We would also like to thank Paul Balmer and Henning Krause for valuable comments on an earlier version of the paper. 
\section{Background on monoidal triangulated categories}
\label{background-mtc}
This section contains background material on monoidal triangulated categories, thick ideals, the noncommutative Balmer spectrum and localization functors. For details the reader is referred to \cite{NVY1, NVY2, NVY3}. 
\subsection{Monoidal triangulated categories}
A {\em{monoidal triangulated category}} ({\em{M$\Delta$C}} for short) is a triangulated category $\bK$ with a monoidal structure consisting of a tensor product $\otimes$ and identity object $\unit$, for which the bifunctor
\[
\otimes : \bK \times \bK \to \bK
\]
is biexact. If the category $\bK$ is $\Bbbk$-linear for a field $\Bbbk$, the tensor product $\otimes$ will be assumed to be $\Bbbk$-bilinear. The associativity and unit constraints of $\bK$ 
will be denoted by
\begin{align*}
&a_{A,B,C} : (A \otimes B) \otimes C \stackrel{\cong}{\rightarrow}
A \otimes (B \otimes C), 
\\
&l_A : \unit \otimes A \stackrel{\cong}{\rightarrow} A
\quad \mbox{and} \quad 
r_A : A \otimes \unit \stackrel{\cong}{\rightarrow} A,
\end{align*}
where $A, B, C \in \bK$. We also have canonical isomorphisms
\[
\Sigma^m (A) \otimes B \cong \Sigma^m(A \otimes B) \cong A \otimes \Sigma^m(B), \quad \forall A, B \in \bK, m \in \mb{Z}
\]
satisfying natural compatibility conditions.

An M$\Delta$C, $\bT$, is called {\em{compactly generated}} if $\bT$ admits arbitrary set indexed coproducts, the tensor product $\otimes$ preserves set indexed coproducts, 
$\bT$ is compactly generated as a triangulated category, the tensor product of compact objects is compact, $\unit$ is a compact object, 
and every compact object is rigid. When $\bT$ is a compactly generated M$\Delta$C, the full subcategory of compact objects of $\bT$ is an M$\Delta$C on its own. It is denoted by $\bT^c$ and is called the compact part of $\bT$. 

\bde{small-big}
An M$\Delta$C, $\bK$, is called {\em{small}} \cite{BS},  if there exists a compactly generated M$\Delta$C, $\widehat{\bK}$, whose
compact part is $\bK$. 
We will say that $\widehat{\bK}$ is an associated {\em{big M$\Delta$C}}
to the small M$\Delta$C, $\bK$. 
\ede
\subsection{Extended endomorphism rings of objects}
The {\em{Tate extended endomorphism ring}} of an object $A$ of $\bK$ is defined to be the set
\begin{equation}
\label{ext-end}
\widehat{\End}^\bullet_{\bK} (A) := \bigoplus_{n \in \mb{Z}} \Hom_\bK(A, \Sigma^n (A)),
\end{equation}
with the product 
\[
f g := \Sigma^n(f) \circ g \in \Hom_\bK(\unit, \Sigma^{m+n} (\unit))
\quad 
\mbox{for}
\quad
f \in \Hom_\bK(A, \Sigma^m (A)), \; 
g \in \Hom_\bK(A, \Sigma^n (A)). 
\]
The {\em{extended endomorphism ring}} $A \in \bK$ is defined to be the subring 
\begin{equation}
\label{ext-endT}
\End^\bullet_{\bK} (A) := \bigoplus_{n \geq 0} \Hom_\bK(A, \Sigma^n (A)).
\end{equation}
When the category $\bK$ is clear from the context, the subscript in the notations in \eqref{ext-end}-\eqref{ext-endT} will be suppressed. 

\subsection{Algebras in $\bK$} It will be useful to introduce the following notion of algebra to formulate the idea of idempotent algebras and orthogonal decomposition in an M$\Delta$C. 

\bde{alg} \hfill 
\begin{enumerate}
\item[(a)] An {\em{algebra}} in $\bK$ is a triple $(A,\mu,\iota)$, where
$A \in \bK$, $\mu \in \Hom_{\bK}(A \otimes A,A)$ and 
$\iota \in \Hom_{\bK}(\unit,A)$ are such that 
\begin{align*}
&\mu \circ (\id_A \otimes \mu) 
\circ a_{A,A,A} = \mu \circ (\mu \otimes \id_A),
\\
&l_A = \mu \circ ( \iota \otimes \id_A) \quad \mbox{and} \quad 
r_A= \mu \circ (\id_A \otimes \iota).
\end{align*}
\item[(b)] An {\em{idempotent algebra}} in $\bK$ is an algebra $(A, \mu, \iota)$ in $\bK$ for which the product map $\mu \in \Hom_{\bK}(A \otimes A, A)$ is an isomorphism. 
\item[(c)] Two algebras $(A_1, \mu_1, \iota_1)$
and $(A_2, \mu_2, \iota_2)$ in $\bK$ are called 
{\em{orthogonal}} if $A_1 \otimes A_2 \cong 0$.
\item[(d)] For an idempotent algebra $(A, \mu, \iota)$, define the 
{\em{corner subcategory}} of $\bK$ generated by $(A, \mu, \iota)$ to be the thick subcategory of $\bK$ generated by the full subcategory with objects $A \otimes C \otimes A$ for $C \in \bK$; the latter is obviously an M$\Delta$C. It is easy to see that the corner subcategory of $\bK$ generated by an idempotent algebra $(A, \mu, \iota)$ is an M$\Delta$C.  
\item[(e)] A {\em{strong idempotent algebra}} in $\bK$ is an idempotent algebra $(A, \mu, \iota)$ for which the objectiwse isomorphisms 
\begin{align*}
&A \otimes ( A \otimes C \otimes A) \stackrel{\mu \otimes \id_C \otimes \id_A}{\longrightarrow} 
A \otimes C \otimes A \quad \mbox{and}
\\
&( A \otimes C \otimes A) \otimes A \stackrel{\id_A \otimes \id_C \otimes \mu}{\longrightarrow} 
A \otimes C \otimes A
\end{align*}
upgrade to unit constraint isomorphisms making $A$ a unit object of the corner subcategory of $\bK$ generated by $(A, \mu, \iota)$.
\end{enumerate}
\ede

By the main result of \cite{SA}, for every strict idempotent algebra 
$(A, \mu, \iota)$, the Tate endomorphism ring of $A$ is commutative, and more precisely, for $f \in \Hom_\bK(A, \Sigma^m (A))$ and $g \in \Hom_\bK(A, \Sigma^n (A))$, the product $f g = \Sigma^n(f) \circ g \in \Hom_\bK(\unit, \Sigma^{m+n} (\unit))$ is also given by 
\[
f g = \mu \circ (f \otimes g) \circ \mu^{-1} = (-1)^{mn} \Sigma^m(g) \circ f. 
\]
More generally, this result holds for any idempotent algebra 
$(A, \mu, \iota)$ in $\bK$ for which $A$ is an identity object of the 
thick triangulated subcategory of $\bK$ generated by $\unit$. 

\subsection{Thick ideals and the noncommutative Balmer spectrum} 
A {\em{thick ideal}} $\bI$ of an M$\Delta$C, $\bK$, is a full triangulated subcategory closed under direct summands such that $\bI$ satisfies the ideal condition: for each $A \in \bI$ 
and $B \in \bK$, the objects $A \otimes B$ and $B \otimes A$ are both in $\bI$. If only $A \otimes B$ (respectively $B \otimes A$) is required to be in $\bI$, then $\bI$ is a {\em{right}} (respectively {\em{left}}) {\em{thick ideal}} of $\bK$. In what follows, if an ideal is only one-sided, it will be explicitly mentioned. For a collection of objects $\bS$ of an M$\Delta$C, $\bK$, denote by $\langle \bS \rangle$ the (unique) minimal thick ideal of $\bK$ containing $\bS$. 

We provide the basic definitions on ideals that will be use throughout the paper. 

\bde{def-ideals} Let $\bK$ be a M$\Delta$C. 
\begin{itemize} 
\item[(a)] Two thick ideals $\bI$ and $\bJ$ of $\bK$ are called {\em coprime} if
\[
\langle \bI, \bJ \rangle = \bK.
\]
\item[(b)] A {\em{principal ideal}} of $\bK$ is a thick ideal generated by a single object, i.e., and ideal of the form $\langle A \rangle$ for some $A \in \bK$. 
\item[(c)] A proper thick ideal $\bP$ of $\bK$ is called {\em{completely prime}} 
if $A \otimes B\in \bP$ implies that one of $A$ or $B$ is an object in $\bP$. 
\item[(d)] A proper thick ideal $\bP$ of $\bK$ is called {\em{prime}} 
if $\bI \otimes \bJ \subseteq \bP$ implies that one of $\bI$ or $\bJ$ is contained in $\bP$, for all thick ideals $\bI$ and $\bJ$ of $\bK$. 
This is equivalent to the property that for all $A, B \in \bK$, 
\[
A \otimes C \otimes B \in \bP, \; \forall C \in \bP \quad \Rightarrow A \in \bP \; \mbox{or} \; B \in \bP,
\]
see \cite[Theorem 1.2.1(a)]{NVY1}.
\item[(e)] A proper thick ideal of $\bK$ is called 
{\em{maximal}} if it is a maximal element of the set of all 
proper thick ideals of $\bK$. All maximal ideals $\bK$ are prime 
\cite[Theorem 3.2.3]{NVY1}.
\end{itemize} 
\ede

{\em{The noncommutative Balmer spectrum of $\bK$}}, denoted by 
$\Spc \bK$, is the collection of all prime ideals of $\bK$, equipped with the Zariski topology where closed sets are defined 
to be arbitrary intersections of the base of closed sets
$$V (A) = \{ \bP \in \Spc \bK : A \not \in \bP \}.$$ 
The map 
\[
V : \bK \to \Spc \bK,
\]
which sends objects of $\bK$ to closed subsets of $\Spc \bK$, is called the {\em{universal support map}}. It is extended to a support map $\Phi$ from the set of thick ideals of $\bK$ to the set of specialization closed subsets of $\Spc \bK$ by setting
\[
\Phi(\bI) := \bigcup_{A \in \bI} V(A).
\]
Every completely prime ideal in an M$\Delta$C is prime. Let $\operatorname{CPSpc} \bK$ be the topological subspace of $\Spc \bK$ consisting of all completely prime ideals of $\bK$.
Its topology is generated by the sets 
\[
V_{CP}(A) = \{ \bP \in \operatorname{CPSpc} \bK \mid A \not \in \bP \}
\]
for $A \in \bK$, 
and one has 
$$ 
V_{CP}(A)=V(A)\cap \operatorname{CPSpc} \bK.
$$ 

In the case of symmetric (or more generally, braided) M$\Delta$Cs, $\bK$, one has $\Spc \bK=\operatorname{CPSpc} \bK$. This is the theory of tensor triangular geometry as first defined and studied by Balmer \cite{Balmer}. 
The nonsymmetric case was developed in \cite{BKS1,NVY1}.

\subsection{The prime radical of an M$\Delta$C} 
A thick ideal $\bP$ of $\bK$ is called {\em{semiprime}} if it is an intersection of a collection of prime ideals of $\bK$. This is equivalent to any of the following two conditions:
\begin{enumerate}
\item[(i)] For all thick ideals $\bI$ of $\bK$, $\bI \otimes \bI \subseteq \bP$ implies that $\bI$ is contained in $\bP$.
\item[(ii)] For all objects $A \in \bK$, 
\[
A \otimes C \otimes A \in \bP, \; \forall C \in \bP \quad \Rightarrow A \in \bP,
\]
\end{enumerate}
see \cite[Theorem 3.4.2]{NVY1}. It now makes sense to consider the smallest semiprime ideal of $\bK$.

\bde{pr-rad} The {\em{prime radical}} of an M$\Delta$C, $\bK$, is the intersection of prime ideals of $\bK$. 
\ede

 Recall that an object $A$ of a monoidal category $\bK$ is {\em{left dualizable}} if there exists an object $A^*$ (called the {\em{left dual of $A$}}), 
together with evaluation and coevaluation maps
\[
\ev: A^* \otimes A \to 1 \quad \mbox{and} \quad
\coev: 1 \to A \otimes A^*,
\] 
such that the compositions
\[
A \xrightarrow{\coev \otimes \id} A \otimes A^* \otimes A \xrightarrow{ \id \otimes \ev} A \quad \mbox{and} \quad
A^* \xrightarrow{\id \otimes \coev} A^* \otimes A \otimes A^* \xrightarrow{\ev \otimes \id} A^*
\]
are the identity maps on $A$ and $A^*$, respectively. The left dual object $A^*$ is unique up to a unique isomorphism, 
\cite[Proposition 2.10.5]{EGNO}. Similarly, one defines the notions of {\em{right dualizable}} objects and their {\em{right duals}}, 
see \cite[Definition 2.10.2]{EGNO}. An object of a monoidal category is {\em{rigid}} if it is both left and right dualizable. One can now make the connection between dualizability and the prime radical. 

\bpr{trivial-prime-rad} \cite[Proposition 4.1.1]{NVY2}
If $\bK$ is a monoidal triangulated category in which every object is either left or right dualizable, then every thick ideal of $\bK$ is semiprime. In particular, the prime radical 
of $\bK$ equals $0$.
\epr

\subsection{Radicals and nilpotent elements} Let $\bf K$ be an M$\Delta$C. Denote by $\bN$ the thick subcategory consisting of nilpotent objects
$$
\{A\in \bK: \ A^{\otimes n}=0\ \text{for some $n \geq 0$}\}.$$ 
If $\bK$ is symmetric (or more generally, braided), then $\bN$ is a thick tensor ideal. 

If $\bK$ admits a completely prime ideal, the {\em completely prime radical}, $\text{CP-Rad} \bK$, is the intersection of all completely prime ideals in $\bK$. 

The following theorem is motivated by classical ring theoretic results, and gives a precise relationship between the prime radical, completely prime radical, and the set of nilpotent elements. 

\bpr{pr-radicals-nilpotents} Let $\bf K$ be an M$\Delta$C, Then one has  the following inclusions: 
\begin{itemize}
\item[(a)] $\operatorname{Rad} {\bK}\subseteq \bN$.
\item[(b)] If $\bK$ admits a completely prime ideal then 
$$\operatorname{Rad} {\bK}\subseteq \bN \subseteq \operatorname{CP-Rad} \bK.$$, 
\item[(c)] If $\bK$ is symmetric then 
$$\operatorname{Rad} {\bK}=\bN=\operatorname{CP-Rad} \bK.$$
\end{itemize} 
\epr

\begin{proof} (a) This is a consequence of \cite[Theorem 3.2.3]{NVY1}. If $A\notin \bN$, take the multiplicative set of all elements of the form $A^{\otimes n}$ for all $n\geq 0$, and take $I=(0)$. Then the theorem guarantees that there is a prime ideal $P$ which does not contain $A$. 

(b) The first inclusion follows from part (a). For the second inclusion, suppose that $A\in \bN$, then $A^{\otimes n}=(0)\in \operatorname{CP-Rad} \bK$ for some $n\geq 0$. Using the properties of completely primeness, 
this implies that $A\in P$ for all $P\in \operatorname{CP-Rad} \bK$, thus $A\in \operatorname{CP-Rad} \bK$.  

(c) This follows from (a) and (b), and the fact that in the symmetric case prime ideals and completely prime ideals coincide. 
\end{proof} 

We should mention that the theorem above encompasses the case of tensor triangular categories. If $\bK$ is a (symmetric) tensor triangular category, it was shown by Balmer that $\operatorname{Rad} \bK=\bN$. See \cite[Corollary 2.4]{Balmer}.

In the case where $\bK$ is rigid and non-commutative, one has $\operatorname{Rad} \bK=(0)$. However, Benson--Witherspoon \cite{Benson-Witherspoon}, and Plavnik--Witherspoon \cite{Plavnik-Witherspoon} have constructed examples of non-cocommutative finite-dimensional Hopf algebras where there are non-trivial nilpotent elements in the stable module category. This demonstrates that there are completely prime ideals that are not prime, and in this situation, the universal support theory does not satisfy the tensor product property. An interesting question would be give necessary and sufficient conditions on when the universal support theory satisfies the tensor product property using the structure theory of M$\Delta$C's involving radicals and nilpotent elements.

\subsection{Localization functors}
\label{sec:loc-funct}

Let $\bK$ be a small M$\Delta$C which is the compact part of $\widehat \bK$, a compactly-generated M$\Delta$C. When $\bI$ is a thick subcategory of $\bK$ (in particular, a thick ideal), we set $\Loc(\bI) \subseteq \widehat \bK$ the localizing subcategory of $\widehat \bK$ generated by $\bI$, that is, the smallest triangulated subcategory of $\widehat \bK$ containing $\bI$ that is closed under arbitrary set-indexed coproducts. By a version of the Eilenberg swindle, $\Loc(\bI)$ is thick, and when $\bI$ is a thick ideal, so is $\Loc(\bI)$. Additionally, we set 
\[
\Loc(\bI)^\perp := \{ A \in \widehat{\bK} : \Hom_{\widehat{\bK}}(B,A) \cong 0 \;\; \forall \;\; B \in \Loc(\bI)\}.
\]

We will denote by $\widehat \bK / \Loc(\bI)$ the Verdier quotient of $\widehat \bK$ by $\Loc(\bI)$; this is a triangulated category whose objects are equal to the objects of $\widehat \bK$, but where the objects isomorphic to 0 are precisely the objects in $\Loc(\bI)$, see \cite[Section 2.1]{Neeman1}. By Brown Representability (c.f.~ \cite[Section 8.2]{Neeman1}), the inclusion functor $\Loc(\bI) \to \widehat{\bK}$ and the projection functor $\widehat{\bK} \to \widehat{\bK}/\Loc(\bI)$ have right adjoints, giving a diagram
\[
\Loc(\bI) \begin{array}{c} {i_*} \atop {\longrightarrow} \\ {\longleftarrow}\atop{i^!} \end{array}   \widehat{\bK} \begin{array}{c} {j^*} \atop {\longrightarrow} \\ {\longleftarrow}\atop{j_*} \end{array}   \widehat{\bK}/ \Loc(\bI).
\]
Using these functors, the following is obtained:
\bth{Rickard-idem}\cite{Rickard1}
If $\widehat \bK$ is a compactly-generated M$\Delta$C with compact part $\bK$ and $\bI$ is a thick subcategory of $\bK$, then there exist functors $\Gamma_{\bI}$ and $L_{\bI}$ from $\widehat{\bK} \to \widehat{\bK}$, which gives for every object $A$ of $\widehat \bK$ a distinguished triangle
\begin{equation}
    \label{local-tri}
\Gamma_{\bI} (A) \to A \to L_{\bI}(A) \to \Sigma (\Gamma_{\bI} (A)),
\end{equation}
such that
\begin{enumerate}
\item[(i)] $\Gamma_{\bI}(A)$ is in $\Loc(\bI)$,
\item[(ii)] $L_{\bI}(A)$ is in $\Loc(\bI)^\perp$, and
\item[(iii)] the distinguished triangle (\ref{local-tri}) is unique up to unique isomorphism.  
\end{enumerate}
\eth
The functors $\Gamma_{\bI}$ and $L_{\bI}$ are referred to as the {\it colocalizing} and {\it localizing} functors associated to $\bI$ respectively, and are formed by the compositions
\[
\Gamma_{\bI} :=i_* \circ i^!, \quad L_{\bI} :=j_* \circ j^*. 
\]
For additional background on localization and colocalization functors, see \cite[Section 3]{BIK2012} and \cite[Section 3]{Stevenson}.

\section{A decomposition theorem for Verdier quotients}
\label{decomposition-thm}
In this section we obtain a Chinese 
Remainder Theorem for the Verdier 
quotients of a big M$\Delta$C, $\widehat{\bK}$, with respect to the localizing ideals of a collection of 
pairwise coprime thick ideals of the corresponding small M$\Delta$C, $\bK$.
\subsection{Statement of the theorem}
For two M$\Delta$Cs $(\bK', \otimes', \unit')$ and $(\bK'', \otimes', \unit'')$, the external direct sum
$\bK' \oplus \bK''$ has a canonical structures of an M$\Delta$C with tensor product
\[
(A' \oplus A'') \otimes (B' \oplus B'') := (A' \otimes B') \oplus (A'' \otimes B'') \quad \mbox{for} \quad A', B' \in \bK', A'', B'' \in \bK''
\]
and unit object $\unit := \unit' \oplus \unit''$.

The following is our Chinese Remainder Theorem for monoidal triangulated categories. Recall the construction of the idempotent completion (Karoubian envelope) of an additive category $\bC$, see e.g. \cite[Section 9.12.1]{EGNO}. Its objects are defined to be pairs $(A,e)$ where $A$ is an object of $\bC$, and $e: A \to A$ is an idempotent morphism. The morphisms between $(A,e)$ and $(A',e')$ are morphisms
$f \in \Hom_\bC(A, A')$ such that $fe = e'f$.
We denote the idempotent completion of $\bC$ by $\bC^{\idem}$. 

The functor $\bC \to \bC^{\idem}$ given by $A \mapsto (A, \id_A)$ embeds $\bC$ as a subcategory of its idempotent completion, which we use without further reference below. If $\bC$ is monoidal, it is clear that $\bC^{\idem}$ is also monoidal, in a way compatible with the embedding $\bC \to \bC^{\idem}$. Furthermore, if $\bT$ is triangulated, then it is a theorem of Balmer and Schlichting that $\bT^{\idem}$ inherits a triangulated structure \cite[Theorem 1.5]{Balmer-Schlichting}. 
\bth{CRT}
Let $\bK$ be a small M$\Delta$C which is the compact part of an associated big 
M$\Delta$C, $\widehat \bK$. 
Assume that $\bI_1, \ldots, \bI_n$ is a collection of thick ideals of $\bK$ that are pairwise coprime, see \deref{def-ideals}(a). 
Then we have the equivalences of M$\Delta$Cs 
\begin{equation}
    \label{CRT-eq1}
\widehat \bK/\Loc(\bI_1 \cap \ldots \cap \bI_n) \cong \widehat \bK / \Loc( \bI_1) \oplus \ldots \oplus \widehat \bK / \Loc(\bI_n)
\end{equation}
and
\begin{equation}
    \label{CRT-eq2}
\left ( \bK/(\bI_1 \cap \ldots \cap \bI_n) \right )^{\idem} \cong (\bK / \bI_1)^{\idem} \oplus \ldots \oplus (\bK / \bI_n)^{\idem}.
\end{equation}
\eth
\subsection{Auxiliary lemmas}
Before proving the theorem, we obtain some auxiliary results. In Lemmas \ref{lLI-rel-tens}--\ref{lcoprime2-3}, we will assume that $\bK$ is a small M$\Delta$C, which is the compact part of an associated big M$\Delta$C, $\widehat \bK$.

\ble{LI-rel-tens}
Let $\bI$ be a thick ideal of $\bK$ and $A,B \in \widehat \bK$. Then 
\[
L_{\bI}(A) \otimes B \cong L_{\bI}(A \otimes B) \cong A \otimes L_{\bI}(B).
\]
\ele

\begin{proof} Consider the distinguished triangle
\[
\Gamma_{\bI} (A) \to A \to L_{\bI}(A) \to \Sigma \Gamma_{\bI}(A).
\]
By exactness of the tensor product, we have a distinguished triangle
\[
\Gamma_{\bI} (A) \otimes B \to A \otimes B \to L_{\bI}(A) \otimes B \to \Sigma \Gamma_{\bI}(A) \otimes B. 
\]
Now $\Gamma_{\bI}(A) \otimes B$ is in the localizing ideal generated by $\bI$, which follows from the fact that localizing categories generated by thick ideals are localizing ideals. On the other hand, we claim that $L_{\bI}(A) \otimes B$ admits no maps from $\Loc(\bI)$ (i.e., $L_{\bI}(A) \otimes B$ is in $\Loc(\bI)^{\perp}$), from which the result will follow by the uniqueness of such a triangle. Note that by duality, this is true if $B$ is compact. Now note that if 
\[B_1 \to B_2 \to B_3 \to \Sigma B_1
\]
is a distinguished triangle such that $L_{\bI}(A) \otimes B_1$ and $L_{\bI}(A) \otimes B_2$ are in $\Loc(\bI)^{\perp}$, then $L_{\bI}(A) \otimes B_3$ is also in $\Loc(\bI)^{\perp}$ using the long exact sequences associated to $\Hom$ functors. It also follows that if $B$ satisfies $L_{\bI}(A) \otimes B$ is in $\Loc(\bI)^{\perp}$, then so is $L_{\bI}(A) \otimes \Sigma B \cong \Sigma(L_{\bI}(A) \otimes B)$. 

Lastly, let $B= \coprod_{i \in I} B_i$ be a set-indexed coproduct. Then if each $B_i$ satisfies $L_{\bI}(A) \otimes B_i \in \Loc(\bI)^{\perp}$, then $B$ also satisfies $L_{\bI}(A) \otimes B \in \Loc(\bI)^{\perp}$, using the fact that the tensor product and $\Hom$ commute with arbtirary coproducts. Therefore, the collection of objects
\[
\{B : \L_{\bI}(A) \otimes B \in \Loc(\bI)^{\perp}\}
\]
is a localizing category. Since it contains all compact objects of $\widehat \bK$, it is in fact equal to $\widehat \bK$. It now follows that $L_{\bI}(A) \otimes B \cong L_{\bI}(A \otimes B)$. The fact that these objects are isomorphic to $A \otimes L_{\bI}(B)$ is analogous.  
\end{proof}
\ble{idem-perp} Let $\bI_1$ and $\bI_2$ be two thick ideals of $\bK$. 
\begin{enumerate}
\item[(a)] For all objects $A, B$ of $\bK$
\[
L_{\bI_1}(A) \otimes L_{\bI_2}(B) \in 
\langle \bI_1, \bI_2 \rangle^{\perp}.
\]
\item[(b)] If $\bI_1$ and $\bI_2$ are coprime
(recall \deref{def-ideals}(a)), then for all objects $A$ and $B$ of $\widehat{\bK}$,
\[
L_{\bI_1}(A) \otimes L_{\bI_2}(B) \cong 0.
\]
\end{enumerate}
\ele

\begin{proof} (a) First, note that for any $A \in \bK$, we have $L_{\bI_1} L_{\bI_2}(A) \in \langle \bI_1, \bI_2 \rangle^\perp$. This follows from the fact that 
    \[
    \Hom_{\widehat{\bK}}(B,L_{\bI_1} L_{\bI_2}(A)) \cong 0
    \]
    for all $B \in \bI_1$ and $\bI_2$. It is clear that it follows that   
    \[
    \Hom_{\widehat{\bK}}(C,L_{\bI_1} L_{\bI_2}(A)) \cong 0
    \]
    for any $C$ in the thick subcategory generated by $\bI_1$ and $\bI_2$, using the long exact sequence of $\Hom$ spaces associated to distinguished triangles. Furthermore, if 
    \[
    \Hom_{\widehat{\bK}}(C,L_{\bI_1} L_{\bI_2}(A)) \cong 0    
    \]
    for any $A \in \bK$, then 
\begin{align*}
    \Hom_{\widehat{\bK}}(C\otimes D,L_{\bI_1} L_{\bI_2}(A)) &\cong \Hom_{\widehat{\bK}}(C,D^{*}\otimes L_{\bI_1} L_{\bI_2} (A) )  \\
    &\cong \Hom_{\widehat{\bK}}(C,L_{\bI_1} L_{\bI_2} (D^{*}\otimes A))\\
    &\cong 0,
\end{align*}
and similar for $D \otimes C$, for arbitrary $D \in \bK$, using \leref{LI-rel-tens}. By another application of \leref{LI-rel-tens}, the result follows.

Part (b) follows at once from the first part. 
\end{proof}

\ble{coprime2-3} If $\bI_1$ and $\bI_2$ are thick ideals of $\bK$, each of which is coprime to a thick ideal $\bI_3$ of $\bK$, then 
$\bI_1 \cap \bI_2$ is coprime to $\bI_3$. 
\ele
\begin{proof} By \cite[Lemma 3.12]{NVY3}, for any two collections of objects $\mc{M}$ and $\mc{N}$ of $\bK$, one has the containment
\[
\langle \mc{M} \rangle \otimes \langle \mc{N} \rangle \subseteq \langle \mc{M} \otimes \bK \otimes \mc{N} \rangle.
\]
Therefore, 
\begin{align*}
\bK \otimes \bK &= \langle \bI_1, \bI_3 \rangle \otimes
\langle \bI_2, \bI_3 \rangle 
\subseteq \langle \bI_1 \otimes \bI_2, \bI_1 \otimes \bI_3, \bI_3 \otimes \bI_2, \bI_3 \otimes \bI_3 \rangle 
\\
&\subseteq \langle \bI_1 \otimes \bI_2, \bI_3 \rangle
\subseteq \langle \bI_1 \cap \bI_2, \bI_3 \rangle. 
\end{align*}
This means that $\unit = \unit \otimes \unit \in \langle \bI_1 \cap \bI_2, \bI_3 \rangle$, and thus 
\[
\langle \bI_1 \cap \bI_2, \bI_3 \rangle = \bK
\]
since the left hand side is a thick ideal of $\bK$. This shows 
that the thick ideals $\bI_1 \cap \bI_2$ and $\bI_3$ are coprime.
\end{proof}
\subsection{Proof of the decomposition theorem}
\hfill \\
{\em{Proof of \thref{CRT}.}} 
We first prove \eqref{CRT-eq1} in the case of two coprime ideals $\bI$ and $\bJ$ of $\bK$.  
Applying $L_{\bJ}$ to the distinguished triangle
\[
\Gamma_{\bI} (A) \to A \to L_{\bI} (A) \to \Sigma \Gamma_{\bI} (A),
\]
and using \leref{idem-perp}(ii), we obtain
\[
L_{\bJ} \Gamma_{\bI} (A) \to L_{\bJ} (A) \to 0 \to \Sigma L_{\bJ} \Gamma_{\bI} (A).
\]
Thus, 
\begin{equation}
\label{L-Gamma}
L_{\bJ} \Gamma_{\bI} (A) \cong L_{\bJ} (A)
\end{equation}
in $\widehat{\bK}$. Now we also have the distinguished triangle
\[
\Gamma_{\bJ} \Gamma_{\bI} (A) \to \Gamma_{\bI} (A) \to L_{\bJ} \Gamma_{\bI} (A)
\]
in $\widehat{\bK}$. Since $\Gamma_{\bJ} \Gamma_{\bI} (A)$ is in $\Loc(\bI) \cap \Loc(\bJ) = \Loc(\bI \cap \bJ)$ (as in \cite[Lemma B.0.2]{NVY3}), we have an isomorphism $\Gamma_{\bI} (A) \cong L_{\bJ} \Gamma_{\bI} (A) \cong L_{\bJ} (A)$ in the Verdier quotient $\widehat \bK / \Loc(\bI \cap \bJ)$. 

Now by working in the Verdier quotient, we have the distinguished triangle
\[
L_{\bJ} (A) \to A \to L_{\bI} (A) \to \Sigma L_{\bJ} (A).
\]
But since $L_{\bI} (A) \cong \Gamma_{\bJ} (A)$, it is in $\Loc(\bJ)$, and $\Sigma L_{\bJ} (A)$ is in $\Loc(\bJ)^{\perp}$, and so the map $L_{\bI} (A) \to \Sigma L_{\bJ} (A)$ is 0. By the Splitting Lemma for triangulated categories, we have an isomorphism of distinguished triangles
\begin{center}
\begin{tikzcd}
L_{\bJ} (A) \arrow[r] \arrow[d, no head, Rightarrow, no head] & A  \arrow[r] \arrow[d, "\cong"]      & L_{\bI} (A) \arrow[r, "0"] \arrow[d, no head, Rightarrow, no head] & \Sigma L_{\bJ} (A) \\
L_{\bJ} (A)  \arrow[r]                               & L_{\bJ} (A) \oplus L_{\bI} (A) \arrow[r] & L_{\bI} (A) \arrow[r]                                     & \Sigma L_{\bJ} (A)
\end{tikzcd}
\end{center}
Hence, in $\widehat \bK  / \Loc(\bI \cap \bJ)$, each object breaks down as a direct sum of objects in $\Loc(\bI)^\perp$ and $\Loc(\bJ)^\perp$. Note that $\Loc(\bI)^\perp$ (resp. $\Loc(\bJ)^{\perp}$) is equivalent to $\widehat \bK / \Loc(\bI)$ (resp. $\widehat \bK/ \Loc(\bJ)$) via $L_{\bI}$ (resp. $L_{\bJ}$). Furthermore, note that $\Loc(\bJ)^{\perp} \otimes \Loc(\bI)^{\perp} = 0$ by \leref{idem-perp}(ii), and there are no nonzero morphisms between objects in $\Loc(\bI)^{\perp}$ and $\Loc(\bJ)^{\perp}$, since for any object $A$, we have (from above) that $\Gamma_{\bI}(A) \cong L_{\bJ} (A)$, and, symmetrically, that $\Gamma_{\bJ} (A) \cong L_{\bI} (A)$. This proves \eqref{CRT-eq1}.

The general case of \eqref{CRT-eq1} follows from the case of $n=2$ by induction and \leref{coprime2-3}. 
 
The isomorphism \eqref{CRT-eq2} follows from the isomorphism 
\eqref{CRT-eq1}, since by \cite[Theorem 5.6.1]{Krause}, if $\bI$ is a thick ideal of $\bK$, then $(\widehat{\bK}/ \Loc(\bI))^c \cong (\bK/ \bI)^{\idem}.$
\qed
\section{Decomposition results for endomorphism algebras}
\label{end-als}
In this section we obtain consequences of 
\thref{CRT} for the decomposition of the 
extended endomorphism algebras 
of idempotent algebra objects of an 
M$\Delta$C and its Verdier quotients. 
\subsection{A decomposition of the extended endomorphism ring of the identity object of a Verdier quotient}
In the setting of \thref{CRT}, 
denote by $\unit'$ and $\unit'_j$ the identity objects of the M$\Delta$C's, 
\[
\bK/\Loc(\bI_1 \cap \ldots \cap \bI_n) \quad \mbox{and} \quad 
\bK/\Loc(\bI_j),
\]
respectively, where $1 \leq j \leq n$. It is easy to see that, under the isomorphism \eqref{CRT-eq2}, we have 
\[
\unit' = \unit'_1 \oplus \ldots \oplus \unit'_n.
\]
Comparing $\Hom$-spaces in \thref{CRT} now implies the following decomposition result. 
\bth{end-algs} In the setting of \thref{CRT}, the decompositions of the extended endomorphism ring and the Tate extended endomorphism ring of the identity object of  
$\widehat \bK/\Loc(\bI_1 \cap \ldots \cap \bI_n)$
in terms of those rings for the identity object of 
$\bK/\Loc(\bI_i)$ are given by:
\[
    \End^\bullet(\unit') \cong \bigoplus_{j=1}^n\End^\bullet(\unit'_j) \quad \mbox{and} \quad
    \widehat{\End}^\bullet(\unit') \cong \bigoplus_{j=1}^n \widehat{\End}^\bullet(\unit'_j),
\]
recall \eqref{ext-end}-\eqref{ext-endT}. 
\eth 
For some related decompositions of $\Hom$-spaces in triangulated categories, see 
\cite[Corollary 5.8]{BalmerFavi0} and \cite[Lemma 3.3]{Krause2}.

Next we investigate the properties of 
special preimages of the unit objects 
$\unit_j$ of $\widehat{\bK}/\Loc(\bI_j)$, 
prove that they give rise to pairwise orthogonal idempotent algebras of $\bK$, show that the unit object of $\bK$ is 
recursively built from them modulo
$\bI_1 \cap \ldots \cap \bI_n$, and describe decomposition results for the corresponding 
(Tate) extended endomorphism rings. 
\subsection{Orthogonal idempotent decomposition of the unit object 
of a small M$\Delta$C}
For a thick ideal $\bI$ of a small M$\Delta$C, $\bK$, 
the isomorphism from \leref{LI-rel-tens} for $A = B:=\unit$ 
gives rise to the isomorphism 
\begin{equation}
    \label{mu-I}
\mu_\bI : L_\bI(\unit) \otimes L_\bI(\unit) 
\stackrel{\cong}{\longrightarrow} 
L_\bI(\unit \otimes L_\bI(\unit) )
\stackrel{L_{\bI}(l_{L_\bI(\unit)})}{\longrightarrow}
L_\bI(L_\bI(\unit) ) \cong L_\bI(\unit).
\end{equation}
\bpr{idemp-alg} Let $\bK$ be a small M$\Delta$C with an associated big M$\Delta$C, $\widehat{\bK}$. 
\begin{enumerate}
\item[(a)] For each thick ideal $\bI$ of $\bK$, 
the triple $(L_\bI(\unit), \mu_\bI, \iota_\bI)$
is a strict idempotent algebra in $\widehat{\bK}$, recall 
\deref{alg}(ii), 
where $\iota_\bI \in \Hom_\bK (\unit, L_\bI(\unit))$ is the 
morphism from the distinguished triangle 
\[
\Gamma_\bI(\unit) \to \unit \stackrel{\iota_\bI}{\longrightarrow}
L_\bI \to \Sigma(\Gamma_\bI).
\]
The corresponding corner subcategory of $\widehat{\bK}$ is $\Loc(\bI)^\perp$.
\item[(b] If $\bI_1$ and $\bI_2$ are coprime thick ideals of $\bK$, then the algebras $(L_{\bI_1}(\unit), \mu_{\bI_1}, \iota_{\bI_1})$
and $(L_{\bI_2}(\unit), \mu_{\bI_2}, \iota_{\bI_2})$
are orthogonal, recall \deref{alg}(c).
\end{enumerate}
\epr
\begin{proof} (a) By the definition of the multiplication map $\mu_{\bI}$, the diagram
\[
\begin{tikzcd}
\Gamma_{\bI} (\unit) \otimes \Gamma_{\bI} (\unit) \arrow[d, "\cong"] \arrow[r, "f \otimes f"] & \unit \otimes \unit \arrow[d, "\cong"] \arrow[r, "g \otimes g"] & L_{\bI} (\unit) \otimes L_{\bI} (\unit) \arrow[d, "\mu_{\bI}"] \\
\Gamma_{\bI} (\unit) \arrow[r, "f"]                                                         & \unit \arrow[r, "g"]                                            & L_{\bI} (\unit)                                             
\end{tikzcd}
\]
commutes, where the triangle
\[
\Gamma_{\bI} (\unit) \xrightarrow{f} \unit \xrightarrow{g} L_{\bI} (\unit) \to \Sigma(\Gamma_{\bI} \unit)
\]
is the distinguished triangle as in (\ref{local-tri}). It is straightforward to verify that the diagrams
\[
\begin{tikzcd}
\Gamma_{\bI} (\unit) \otimes \Gamma_{\bI} (\unit) \otimes \Gamma_{\bI} (\unit) \arrow[d, "\cong"] \arrow[r, "f \otimes f \otimes f"] & \unit \otimes \unit \otimes \unit \arrow[d, "\cong"] \arrow[r, "g \otimes g \otimes g"] & L_{\bI} (\unit) \otimes L_{\bI} (\unit) \otimes L_{\bI} (\unit) \arrow[d, "\mu_{\bI} \circ (\mu_{\bI} \otimes \id_{L_{\bI} \unit})"] \\
\Gamma_{\bI} (\unit) \arrow[r, "f"]                                                                                              & \unit \arrow[r, "g"]                                                                    & L_{\bI} (\unit)                                                                                                                 
\end{tikzcd}
\]
and
\[
\begin{tikzcd}
\Gamma_{\bI} (\unit) \otimes \Gamma_{\bI} (\unit) \otimes \Gamma_{\bI} (\unit) \arrow[d, "\cong"] \arrow[r, "f \otimes f \otimes f"] & \unit \otimes \unit \otimes \unit \arrow[d, "\cong"] \arrow[r, "g \otimes g \otimes g"] & L_{\bI} (\unit) \otimes L_{\bI} (\unit) \otimes L_{\bI} (\unit) \arrow[d, "\mu_{\bI} \circ ( \id_{L_{\bI} \unit} \otimes \mu_{\bI})"] \\
\Gamma_{\bI} (\unit) \arrow[r, "f"]                                                                                              & \unit \arrow[r, "g"]                                                                    & L_{\bI} (\unit)                                                                                                                  
\end{tikzcd}
\]
commute. Then by \thref{Rickard-idem}(c), $\mu_{\bI} \circ (\id_{L_{\bI} \unit} \otimes \mu_{\bI}) = \mu_{\bI} \circ (\mu_{\bI} \otimes \id_{L_{\bI} \unit}).$ Likewise, the diagram
\[
\begin{tikzcd}
\unit \otimes \Gamma_{\bI} (\unit) \arrow[d, "\cong"] \arrow[dd, "l_{\Gamma_{\bI} (\unit)}"', bend right=70] \arrow[r, "\id_{\unit} \otimes f"] & \unit \otimes \unit \arrow[d, "\cong"'] \arrow[r, "\id_{\unit} \otimes g"] & \unit \otimes L_{\bI} (\unit) \arrow[d, "\iota_{\bI} \otimes \id_{L_{\bI} (\unit)}"'] \arrow[dd, "l_{L_{\bI} (\unit)}", bend left=70] \\
\Gamma_{\bI} (\unit) \otimes \Gamma_{\bI} (\unit) \arrow[d, "\cong"] \arrow[r, "f \otimes f"]                                                   & \unit \otimes \unit \arrow[d, "\cong"'] \arrow[r, "g \otimes g"]           & L_{\bI} (\unit) \otimes L_{\bI} (\unit) \arrow[d, "\mu_{\bI}"']                                                                     \\
\Gamma_{\bI} (\unit) \arrow[r, "f"]                                                                                                           & \unit \arrow[r, "g"]                                                       & L_{\bI} (\unit)                                                                                                                  
\end{tikzcd}
\]
commutes, and again by \thref{Rickard-idem}(c), $l_{L_{\bI} (\unit)} = \mu_{\bI} \circ (\iota_{\bI} \otimes \id_{L_{\bI} (\unit)}).$ The identity $r_{L_{\bI} (\unit)} = \mu_{\bI} \circ (\id_{L_{\bI} (\unit)} \otimes \iota_{\bI})$ is similar. This shows that 
$(L_\bI(\unit), \mu_\bI, \iota_\bI)$ is an idempotent algebra in $\widehat{\bK}$.

The corner subcategory of $\widehat{\bK}$ generated by $L_\bI(\unit)$ is 
$\Loc(\bI)^\perp$ because the thick triangulated subcategory of $\widehat{\bK}$ 
generated by 
\[
L_\bI(\unit) \otimes C \otimes L_\bI (\unit) \cong L_\bI (C)
\]
for $C \in \widehat{\bK}$ is $\Loc(\bI)^\perp$. The functorial isomorphisms from \leref{LI-rel-tens} show that $L_\bI(\unit)$ is a unit object of $\Loc(\bI)^\perp$, so 
$(L_\bI(\unit), \mu_\bI, \iota_\bI)$ is a strict idempotent algebra in $\widehat{\bK}$.

Part(b) follows from \leref{idem-perp}(b).
\end{proof}

If $\{\bI_1, \ldots, \bI_n\}$
is a collection of pairwise coprime thick ideals of a small M$\Delta$C, $\bK$, 
with a corresponding big 
M$\Delta$C, $\widehat{\bK}$, then we obtain 
the collection of pairwise orthogonal idempotent algebras 
\[
(L_{\bI_j}(\unit), \mu_{\bI_j}, \iota_{\bI_j}), \quad 1 \leq j \leq n 
\]
in $\widehat{\bK}$. 
The next result shows that the unit object of $\widehat{\bK}$ can be recursively built (in any order) from the underlying objects of these algebras modulo the localizing ideal 
\[
\Loc(\bI_1 \cap \ldots \cap \bI_n).
\]
\bpr{recursive} In the setting of \thref{CRT},
for all $\sigma \in S_n$ and $1 \leq k \leq n$,
\begin{multline*}
   \Gamma_{\bI_{\sigma(k)}} \Gamma_{\bI_{\sigma(k-1)}} \ldots 
   \Gamma_{\bI_{\sigma(1)}} (\unit) \to 
   \Gamma_{\bI_{\sigma(k-1)}} \ldots 
   \Gamma_{\bI_{\sigma(1)}} (\unit) 
\to L_{\bI_{\sigma(k)}} (\unit) \\
\to \Sigma \Gamma_{\bI_{\sigma(k)}} \Gamma_{\bI_{\sigma(k-1)}} \ldots 
   \Gamma_{\bI_{\sigma(1)}} (\unit)
\end{multline*}
is a distinguished triangle.
\epr
\begin{proof} Equation \eqref{L-Gamma} implies that 
\[
L_{\bI_{\sigma(k)}} \Gamma_{\bI_{\sigma(k-1)}} \ldots 
   \Gamma_{\bI_{\sigma(1)}} (\unit)
\cong L_{\bI_{\sigma(k)}} (\unit). 
\]
The proposition now follows from this isomorphism by using the localizing distinguished triangles. 
\end{proof}
\subsection{Hom spaces in Verdier quotients} The following result demonstrates how the localization functors behave with respect of the Hom-functors. 

\bpr{hom-mod} Let $\bI$ be a thick ideal of a 
small M$\Delta$C, $\bK$, with an associated big 
M$\Delta$C, $\widehat{\bK}$. 
\begin{itemize}
\item[(a)] Then for all objects $A, B \in \widehat \bK$,
\[\Hom_{\widehat{\bK}/ \Loc(\bI)} (A,B) \cong \Hom_{\widehat{\bK}}(A, L_{\bI}(B)) \cong \Hom_{\widehat{\bK}} (L_\bI(A), L_{\bI}(B)).
\]
\item[(b)] One has the following isomorphisms: 
\[
\End^\bullet_{\widehat{\bK}/ \Loc(\bI)}(\unit) 
\cong \End^\bullet_{\widehat{\bK}} (L_{\bI}(\unit) ) \quad 
\mbox{and} \quad 
\widehat{\End}^\bullet_{\widehat{\bK}/ \Loc(\bI)}(\unit) 
\cong \widehat{\End}^\bullet_{\widehat{\bK}} (L_{\bI}(\unit)).
\]
\end{itemize}
\epr
\begin{proof} (a) Denote for brevity 
\[
\widehat \bM := \widehat{\bK}/ \Loc(\bI). 
\]
First, note that we have a long exact sequence
\[... \to \Hom^i_{\widehat{\bM}}(A, \Gamma_\bI (B)) \to \Hom^i_{\widehat{\bM}}(A, B) \to \Hom^i_{\widehat{\bM}}(A, L_{\bI}(B)) \to \Hom^{i+1}_{\widehat{\bM}}(A, \Gamma_{\bI}(B)) \to ...\]
Since $\Sigma^i \Gamma_{\bI}(B) \cong 0$ in $\widehat \bM$ for all $i$, we have $\Hom_{\widehat{\bM}}(A, L_{\bI}(B)) \cong \Hom_{\widehat{\bM}}(A, B)$. The first isomorphism now follows from \cite[Lemma 4.8.1(3)]{Krause}. The second isomorphism follows by applying the long exact sequence associated to $\Hom_{\widehat \bK}(-, L_{\bI} (B))$ to the distinguished triangle
\[
\Gamma_\bI A \to A \to L_{\bI} A,
\]
and using the fact that $\Hom^i_{\widehat \bK} (\Gamma_\bI( A), L_\bI (B)) =0$ for all $i \in \mathbb{Z}$. 

Part (b) is a direct consequence of (a). 
\end{proof}

\subsection{A decomposition of (Tate) extended endomorphism rings of idempotent algebra objects}
Our next result describes a decomposition of the (Tate) extended endomorphism rings of the idempotent algebra objects 
of $\widehat{\bK}$ that are preimages of the 
unit objects of the Verdier quotients of the 
form 
\[
\widehat{\bK}/ \Loc(\bI_1 \cap \ldots \cap \bI_n).
\]

The following theorem follows from \thref{end-algs} and \prref{hom-mod}(a).  

\bth{decomp-end-idemp-alg} In the setting of \thref{CRT}, we have the following decompositions of the extended endomorphism ring and the Tate extended endomorphism rings of the idempotent algebra object 
$L_{\bI_1 \cap \ldots \cap \bI_n}(\unit)$
of $\widehat{\bK}$:
\[
\End_{\widehat{\bK}}^\bullet
\big(L_{\bI_1 \cap \ldots \cap \bI_n}(\unit) \big) \cong \bigoplus_{j=1}^n \End_{\widehat{\bK}}^\bullet
\big(L_{\bI_j}(\unit) \big), 
\quad
\widehat{\End}_{\widehat{\bK}}^\bullet
\big( L_{\bI_1 \cap \ldots \cap \bI_n} (\unit) \big) \cong
\bigoplus_{j=1}^n 
\widehat{\End}_{\widehat{\bK}}^\bullet
\big(L_{\bI_j}(\unit) \big). 
\]
For the statement, recall \eqref{ext-end}-\eqref{ext-endT}.
\eth

\subsection{Internal (Tate) extended endomorphism algebras}
It is natural to ask whether the (Tate) extended endomorphism rings in \thref{end-algs} have an internal presentation. In \prref{internal-end} we show that the answer to this question is affirmative. 
\bde{internal-hom} Let $\bI$ be a thick ideal of a 
small M$\Delta$C, $\bK$, with an associated big M$\Delta$C, $\widehat{\bK}$. We say that $A, B \in \widehat{\bM}:= \bK/ \Loc(\bI)$ have an internal hom in $\bK$ if there is an object 
\[
\underline{\Hom}_{\widehat{\bM}}(A,B) \in \widehat{\bK} 
\]
and for $C \in \widehat{\bK}$ there exists a natural isomorphism  
\[
\Hom_{\widehat{\bM}}(C \otimes A, B) \cong \Hom_{\widehat{\bK}}(C, \underline{\Hom}_{\widehat{\bM}}(A,B)).
\]
\ede
\bpr{internal-end} In the setting of \deref{internal-hom}, the following hold:
\begin{enumerate}
\item[(a)] The images of every $A \in \bK$ and $B \in \widehat{\bK}$ in $\widehat{\bM}$ have an internal hom in $\widehat{\bK}$, and it is given by 
\[
\underline{\Hom}_{\widehat{\bM}}(A,B) := L_\bI(A^*) \otimes B. 
\]
\item[(b)] The object
\[
\bigoplus_{n \in \mb{Z}} \Sigma^n (L_\bI(\unit)) 
\]
is an algebra in $\widehat{\bK}$ with product induced by the morphisms
\[
\Sigma^m (L_\bI(\unit))  \otimes \Sigma^n (L_\bI(\unit))
\cong \Sigma^{m+n} (L_\bI(\unit) \otimes L_\bI(\unit)) \stackrel{\Sigma^{m+n}(\mu_\bI)}{\longrightarrow}
\Sigma^{m+n} (L_\bI(\unit))
\]
for $m, n \in \mb{Z}$, recall \eqref{mu-I}, and unit 
\[
\unit \stackrel{\iota_{\bI}}{\longrightarrow} L_\bI(\unit) \to \bigoplus_{n \in \mb{Z}} \Sigma^n (L_\bI(\unit)).
\]
\item[(c)] The object 
\[
\bigoplus_{n \geq 0} \Sigma^n (L_\bI(\unit)) 
\]
is an algebra in $\widehat{\bK}$ with similarly defined product and unit. 
\end{enumerate}
\epr
\begin{proof} (a) We must show that $L_{\bI}(A^*) \otimes B$ satisfies the property that 
\[
\Hom_{\widehat{\bM}}(C \otimes A, B) \cong \Hom_{\widehat{\bK}}(C, L_{\bI}(A^*) \otimes B).
\]
This follows from \leref{LI-rel-tens} and \prref{hom-mod}(a). 

(b) The statement follows from \prref{idemp-alg} and from the more general assertion that if $(A, \mu, \iota)$ is an algebra, then there is an induced algebra structure on 
\[
\widehat A := \bigoplus_{i \in \mathbb{Z}} \Sigma^i (A),
\]
where the multiplication is induced via
\[
\Sigma^i (A) \otimes \Sigma^j (A) \xrightarrow{\cong} \Sigma^{i+j} (A) \xrightarrow{ \Sigma^{i+j} \mu} \Sigma^{i+j} (A),
\]
and the unit is induced via
\[
\unit \xrightarrow{\iota} A \hookrightarrow \widehat A.
\]
Associativity and the unit axiom of the multiplication of $\widehat A$ follow from the naturality of the isomorphism $\Sigma^i A \otimes \Sigma^j A \cong \Sigma^{i+j} A$.

Part (c) follows from the statement of (b).     
\end{proof}

\section{Existence of pairs of coprime proper thick ideals}
In this section we give a characterization of when an M$\Delta$C, $\bK$, contains a pair of coprime proper 
thick ideals in terms of the topological properties of the noncommutative Balmer spectrum of $\bK$. 

For this section, we will not impose the standard assumptions. For example, it will not be assumed that $\bK$ has a generator, that it is rigid, that it is the compact part of a compactly generated monoidal triangulated category, or that $\Spc \bK$ is Noetherian. In particular, this means that we do not assume Thomason subsets of the Balmer spectrum of $\bK$ parametrize the thick ideals of $\bK$ (cf.~ \prref{bij}). 

\subsection{Auxiliary lemmas} The following lemmas will be later used to prove the main result in this section (i.e., \thref{coprime-equiv}). 
\ble{closure} For every prime ideal $\bP$ of an M$\Delta$C, $\bK$, the closure of $\bP$ in $\Spc \bK$ is equal to the set of all thick primes ideals $\bQ$ of $\bK$ contained in $\bP$.
\ele

\begin{proof}
    On one hand, if $\bQ \subseteq \bP$ and $\bP \in V(\mc{S})$ for some collection $\mc{S}$ of objects of $\bK$, then it is clear that $\bQ \in V(\mc{S})$ as well, since $V(\mc{S})$ consists of primes which intersect trivially with $\mc{S}$. On the other hand, if $\bQ \not \subseteq \bP$, then $\bP \in V(A)$ and $\bQ \not \in V(A)$, where $A$ is any object in $\bQ$ and not in $\bP$. 
\end{proof}

\ble{proper-supp-proper-id}
A thick ideal $\bI$ of an M$\Delta$C, $\bK$, is proper if and only if $\Phi(\bI)$ is proper in $\Spc \bK$. 
\ele

\begin{proof}
If $\Phi(\bI)= \Spc \bK$, then by the definition of universal support, there are no prime ideals of $\bK$ containing $\bI$. Since maximal ideals of $\bK$ are prime, this condition is satisfied if and only if $\bI = \bK$. 
\end{proof}

\ble{coprime-ideal-support}
Two thick ideals $\bI$ and $\bJ$ of an M$\Delta$C, $\bK$, 
are coprime if and only if 
$\Phi(\bI) \cup \Phi(\bJ) = \Spc \bK$.
\ele
\begin{proof}
    Recall that $\Phi(\langle \bI \cup \bJ \rangle)= \Phi(\bI) \cup \Phi(\bJ)$. Hence, if $\bI$ and $\bJ$ are coprime, then it is clear that $\Phi(\bK)=\Phi(\langle \bI \cup \bJ \rangle)= \Phi(\bI) \cup \Phi(\bJ) = \Spc \bK$. On the other hand, if $\Phi(\langle \bI \cup \bJ \rangle) = \Spc \bK$, then by \leref{proper-supp-proper-id}, $\bI$ and $\bJ$ are coprime.  
\end{proof}
\subsection{} In the following theorem, we prove a characterization theorem for existence of pairs of coprime proper thick ideals. 

\bth{coprime-equiv} For every M$\Delta$C, $\bK$, the following are equivalent:
\begin{itemize}
    \item[(a)] There exist two coprime proper thick ideals of $\bK$.
    \item[(b)] There exist two distinct maximal ideals of $\bK$.
    \item[(c)] There exist two coprime proper principal ideals of $\bK$.
    \item[(d)] $\Spc \bK$ is reducible, that is, there exist proper closed sets $Z_1$ and $Z_2$ in $\Spc \bK$ with $\Spc \bK=Z_1 \cup Z_2.$
\end{itemize}
\eth
\begin{proof}
It is clear that (b) $\Rightarrow$ (a) and (c) $\Rightarrow$ (a). 

(a) $\Rightarrow$ (c): if $\bI$ and $\bJ$ are proper coprime ideals of $\bK$, then $\unit$ can be formed from elements of $\bI$ and $\bJ$ in finitely many steps by taking cones, shifts, direct summands, and tensor products with arbitrary elements of $\bK$. In particular, only finitely many objects of $\bI$ and $\bJ$ need to be used to generate $\unit$; by taking the direct sum of these finite objects from $\bI$ and $\bJ$, respectively, we produce $A \in \bI$ and $B \in \bJ$ so that $\langle A, B \rangle = \bK$. Hence $\langle A \rangle$ and $\langle B \rangle$ are a pair of principal proper coprime ideals.

(c) $\Rightarrow$ (d): by \leref{coprime-ideal-support} and the fact that the universal support of a single object is closed, if $\langle A \rangle$ and $\langle B \rangle$ are proper coprime principal ideals, then $V(A) \cup V(B)$ forms a decomposition of $\Spc \bK$ into closed subsets.

(d) $\Rightarrow$ (b): suppose $\bK$ has a unique maximal ideal $\bM$. If $\Spc \bK= Z_1 \cup Z_2$ for some closed subsets $Z_1$ and $Z_2$, it follows that $\bM \in Z_1$, without loss of generality. But the closure of $\bM$ in $\Spc \bK$ consists of all primes contained in $\bM$, by \leref{closure}, which, since $\bM$ is the unique maximal ideal, consists of $\Spc \bK$. Hence, $Z_1 = \Spc \bK$, and is not proper. 
\end{proof}

As mentioned above, \thref{coprime-equiv} does not use any bijection between thick ideals of $\bK$ and Thomason closed sets of $\Spc \bK$, since (without additional assumptions) such a bijection might not exist. In particular, if $\Spc \bK= Z_1 \cup Z_2$ for closed sets $Z_1$ and $Z_2$, it is not necessarily true that the coprime proper principal ideals $\bI_1$ and $\bI_2$, guaranteed to exist by \thref{coprime-equiv}(c), satisfy $\Phi(\bI_1)=Z_1$ and $\Phi(\bI_2)=Z_2$. The closed sets $\Phi(\bI_1)$ and $\Phi(\bI_2)$ are a pair of proper closed subsets whose union is all of $\Spc \bK$, but they might not equal $Z_1$ and $Z_2$. 

\section{Pairs of complementary comprime thick ideals}
In this section we obtain a topological characterization of 
when an M$\Delta$C, $\bK$, contains a pair of coprime proper thick ideals which are complementary in the sense that their 
intersection lies in the prime radical of $\bK$. This is in turn used to obtain a topological characterization of when a small M$\Delta$C has a direct sum decomposition (as in \thref{CRT}). 
\subsection{Characterization theorem for existence of pairs of coprime complementary thick ideals}

\bde{complem}
We say that two thick ideals of an M$\Delta$C, $\bK$, are {\it complementary} if their intersection is contained in the prime radical of $\bK$. 
\ede

It follows from \prref{trivial-prime-rad} that, if every object in $\bK$ is either left or right dualizable, then 
two ideals of $\bK$ are complementary if and only if their intersection is $\{0\}$. This in particular holds if $\bK$ 
is rigid.

\bth{coprime-complem-equiv}
For every M$\Delta$C, $\bK$, the following are equivalent:
\begin{enumerate}
    \item[(a)] There exist two proper coprime complementary thick ideals of $\bK$.
    \item[(b)] There exist two proper coprime complementary principal ideals of $\bK$.
\end{enumerate}
These conditions imply 
\begin{enumerate}
    \item[(c)]$\Spc \bK$ is disconnected, that is, there exist proper disjoint closed sets $Z_1$ and $Z_2$ in $\Spc \bK$ with $\Spc \bK=Z_1 \sqcup Z_2.$
\end{enumerate}
If either of the two conditions 
\begin{enumerate}
    \item[(i)] $\bK$ is generated by an object $G$ as a thick subcategory, or
    \item[(ii)] $\Spc \bK$ is Noetherian
\end{enumerate}
are satisfied, then conditions (a), (b) and (c) are equivalent.
\eth

\begin{proof}
(b) $\Rightarrow$ (a) is clear. 

(a) $\Rightarrow$ (b): If $\bI$ and $\bJ$ are complementary coprime proper thick ideals, then just as in the proof of \thref{coprime-equiv}, there exist $A \in \bI$ and $B \in \bJ$ with $\langle A, B \rangle = \bK$. Since $\bI \cap \bJ$ is conained in the prime radical, we also have $\langle A \rangle \cap \langle B \rangle$ contained in the prime radical. 

(b) $\Rightarrow$ (c): If $\langle A \rangle$ and $\langle B \rangle$ are coprime complementary proper principal ideals, then
    \[
    V(A) \cap V(B) = \Phi(\langle A \rangle \cap \langle B \rangle)  = \varnothing,
    \] by the assumption that $\langle A \rangle \cap \langle B \rangle$ is in the prime radical of $\bK$. Hence in this case, $\Spc \bK$ is disconnected via
    \[
    V(A) \sqcup V(B).
    \]

    (c) $\Rightarrow$ (b) when $\bK$ satisfies either (i) or (ii): if $\Spc \bK$ is Noetherian, then by \cite[Corollary 5.5]{Rowe}, we have $Z_1=V(A)$ and $Z_2=V(B)$ for some objects $A$ and $B$ in $\bK$. On the other hand, $\bK$ has a generator $G$, then $\bK$ is quasicompact by \leref{qc-open}, since $\Spc \bK = V(\unit)$. It follows that $Z_1$ and $Z_2$ must be quasicompact as well, and again by \leref{qc-open}, $Z_1=V(A)$ and $Z_2=V(B)$ for some objects $A$ and $B$ in $\bK$. In either case, we have realized $Z_1$ and $Z_2$ as the universal supports of objects $A$ and $B$, respectively. Then since $\Phi(\langle A \rangle \cap \langle B \rangle)=V(A) \cap V(B)=Z_1 \cap Z_2= \varnothing$, it follows that $\langle A \rangle$ and $\langle B \rangle$ are complementary.
\end{proof}

In the commutative case, see \cite[Theorem A.5]{BIK2013} for a related decomposition theorem.

\subsection{Interpretation of the decomposition theorem on the level of Balmer spectra}
\thref{CRT} raises the question of when an M$\Delta$C has a direct sum decomposition. Next we give a characterization of this on the level of Balmer spectra. To do this, we first check the following noncommutative analogue of \cite[Proposition 3.13]{Balmer}. Recall the definition of the idempotent completion $\bK^{\idem}$ of $\bK$, as in Section \ref{decomposition-thm}.
\ble{idem-spc}
If $\bK$ is a monoidal triangulated category, then we have $\Spc \bK \cong \Spc \bK^{\idem}$.
\ele

\begin{proof}
The proof follows in a similar manner as the commutative case; we include the details of the noncommutative generalization for the ease of the reader. Recall that by the construction of $\bK^{\idem}$, we have $\bK$ embedded as a full monoidal-triangulated subcategory of $\bK^{\idem}$. Furthermore, for any object $A$ of $\bK^{\idem}$, there exists an object $B$ so that $A \oplus B \in \bK$. We claim that the map
\begin{align*}
    \Spc \bK^{\idem} & \xrightarrow{i} \Spc \bK\\
    \bP &\mapsto \bP \cap \bK
\end{align*}
is a homeomorphism. 

The fact that $i$ is well-defined follows immediately from the fact that $\bK$ generates $\bK^{\idem}$ as a thick subcategory. In more detail, if $\bP$ is in $\Spc \bK^{\idem}$, and $A \otimes \bK \otimes B \subseteq \bP \cap \bK$ for some objects $A$ and $B$ of $\bK$, then $A \otimes \bK^{\idem} \otimes B \subseteq \bP$, and so one of $A$ or $B$ is in $\bP \cap \bK$. The facts that $i(\bP)$ is a thick ideal and that $i$ is continuous are straightforward.

Next, note that if $A$ is in $\bK^{\idem}$, then $A \oplus \Sigma A \in \bK$ by the following observation: there exists $B$ with $A \oplus B \in \bK$, and by adding standard distinguished triangles we obtain a distinguished triangle
\[
A \oplus B \to A \oplus \Sigma (A) \to \Sigma(A \oplus B) \to \Sigma( A \oplus B).
\]
Since the first and third objects of this triangle are in $\bK$, so is the second. 
\vskip .25cm 
\noindent
{\em Injectivity of $i$}: if $\bP \cap \bK = \bQ \cap \bK$ for some prime ideals $\bP$ and $\bQ$ of $\bK^{\idem}$, and if $A$ is in $\bP$, then $A \oplus \Sigma A \in \bP \cap \bK = \bQ \cap \bK$, and so $A \in \bQ$. 
\vskip .25cm 
\noindent
{\em Surjectivity of $i$}: let $\bP \in \Spc \bK$. Set 
\[
\bQ :=\{A \in \bK^{\idem} : A \oplus \Sigma (A) \in \bP\}.
\]
Clearly, $\bQ \cap \bK = \bP$, so we are left with showing that $\bQ$ is in $\Spc \bK^{\idem}$. It is immediate that $\bQ$ is a thick subcategory. To check that $\bQ$ is an ideal, suppose $A$ satisfies $A \oplus \Sigma A \in \bP$. Then note that for any $B \in \bK^{\idem}$, we have 
\[
(A \oplus \Sigma (A)) \otimes (B \oplus \Sigma (B)) = (A \otimes B) \oplus \Sigma(A \otimes B) \oplus \Sigma(A \otimes B) \oplus \Sigma^2 (A \otimes B).
\]
The left hand side is an object of $\bP$ by assumption, so $(A \otimes B) \oplus \Sigma(A \otimes B)$ is as well, as a summand of the right hand side. Thus $A \otimes B \in \bQ$, and so $\bQ$ is a right ideal; the left ideal property follows similarly.
Now, suppose $A \otimes \bK^{\idem} \otimes B \subseteq \bQ$ for some objects $A$ and $B$ of $\bK^{\idem}$. For any $C \in \bK^{\idem}$, 
\[
(A \otimes C \otimes B) \oplus \Sigma(A \otimes C \otimes B) \in \bP.
\]
If $C \in \bK$, then note that this implies that 
\begin{align*}
&(A \oplus \Sigma (A)) \otimes C \otimes (B \oplus \Sigma (B)) 
\\
&= (A \otimes C \otimes B) \oplus \Sigma(A \otimes C \otimes B) \oplus \Sigma(A \otimes C \otimes B) \oplus \Sigma^2 (A \otimes C \otimes B)\\
&=(A \otimes C \otimes B) \oplus \Sigma (A \otimes C \otimes B) \oplus \Sigma \left ( (A \otimes C \otimes B) \oplus \Sigma (A \otimes C \otimes B) \right)
\end{align*}
is in $\bP$. By primeness of $\bP$, we have either $A \oplus \Sigma (A)$ or $B \oplus \Sigma (B) \in \bP$, and hence either $A$ or $B$ is in $\bQ$.
\end{proof}

\bth{dis-spc-sum-cat}
Suppose that $\bK$ is a small M$\Delta$C which is the compact part of an associated big M$\Delta$C, $\widehat \bK$, and that $\bK$ satisfies either
\begin{enumerate}
    \item[(a)] $\bK$ is generated by an object $G$ as a thick subcategory, or
    \item[(b)] $\Spc \bK$ is Noetherian.
\end{enumerate}
If $\Spc \bK$ is disconnected via $\Spc \bK = Z_1 \sqcup Z_2$, then $\widehat \bK$ decomposes as a direct sum $\widehat \bK_1 \oplus \widehat \bK_2$, with the conditions that $\bK_i:=\widehat \bK_i^c$ and $\Spc (\bK_i) \cong Z_i$.  
\eth
\begin{proof}
By Theorems \ref{tCRT} and \ref{tcoprime-complem-equiv}, $\widehat \bK \cong \widehat \bK/\bL_1 \oplus \widehat \bK/\bL_2$ for two thick ideals $\bI_1,$ $\bI_2$ with $\Phi(\bI_i)=Z_i$ and setting $\bL_i:=\Loc(\bI_i)$. Set $\widehat \bK_i :=\widehat \bK / \bL_i$. We have $\widehat \bK_i^c \cong (\bK/\bI_i)^{\idem}$. By \leref{idem-spc}, the Balmer spectrum is not affected by taking the idempotent completion. To complete the proof, it is straightforward to verify that for any thick ideal $\bI$, there is a homeomorphism
\begin{align*}
\Spc (\bK/\bI) & \cong \{ \bP \in \Spc \bK : \bI \subseteq \bP\} = (\Spc \bK) \backslash \Phi(\bI)\\
\bQ & \mapsto \{A \in \bK : q(A) \in \bQ\}\\
\{q(A) : A \in \bP \} & \mapsfrom \bP
\end{align*}
where $q$ is the quotient functor $\bK \to \bK/\bI$. Hence, in the situation of this lemma, we see that $\Spc (\bK/\bI_1) \cong Z_2$ and $\Spc (\bK/\bI_2) \cong Z_1.$\end{proof}



\section{A General Version of Carlson's Theorem}
\label{Carlson}

Our goal is to derive a generalized Carlson theorem about connectedness of the universal support of an indecomposable object of a small M$\Delta$C.

\subsection{An auxiliary decomposition result} Carlson proved in \cite{Carlson1984} the celebrated theorem that the support variety associated to an indecomposable module for a finite group in characteristic $p$ is connected. For symmetric monoidal tensor categories, a version of Carlson's theorem was proven by Balmer \cite[Main Theorem]{Balmer2} and a more general version on smashing ideals was obtained by Balmer--Favi \cite[Corollary 3.15]{BalmerFavi}. In order to state our generalization, we first need another decomposition lemma. 

\ble{decomp-quot} Let $\bK$ be a small M$\Delta$C which is the compact part of an associated big M$\Delta$C, $\widehat \bK$, and $\bI_1$ and $\bI_2$ be two thick ideals of $\bK$. 
Set $\bI:= \langle \bI_1, \bI_2 \rangle $. If $A$ is in $\bI$, then $A \cong L_{\bI_1}(A) \oplus L_{\bI_2}(A)$ in the quotient $\widehat{\bK}/\Loc(\bI_1 \cap \bI_2)$. 
\ele

\begin{proof}
By the distinguished triangle
\[
\Gamma_{\bI_1} (A) \to A \to L_{\bI_1} (A) \to \Sigma \Gamma_{\bI_1} (A),
\]
and the fact that $A \in \bI$ and $\bI_1 \subseteq \bI$, we have $L_{\bI_1}A \in \Loc(\bI)$. Applying $L_{\bI_2}$ to this distinguished triangle,
we obtain
\[
L_{\bI_2} \Gamma_{\bI_1} (A) \to L_{\bI_2} (A) \to 0 \to \Sigma L_{\bI_2} \Gamma_{\bI_1} (A)
\]
by \leref{idem-perp}(a) and the fact that $L_{\bI_2} (A) \in \Loc(\bI)$. It follows $L_{\bI_2} \Gamma_{\bI_1} (A) \cong L_{\bI_2} (A)$ in $\widehat{\bK}$. Now we also have the distinguished triangle
\[
\Gamma_{\bI_2} \Gamma_{\bI_1} (A) \to \Gamma_{\bI_1} (A) \to L_{\bI_2} \Gamma_{\bI_1} (A)
\]
in $\widehat{\bK}$. Since $\Gamma_{\bI_2} \Gamma_{\bI_1} (A)$ is in $\Loc(\bI_1) \cap \Loc(\bI_2) = \Loc(\bI_1 \cap \bI_2)$ (as in \cite[Lemma B.0.2]{NVY3}), we have an isomorphism $\Gamma_{\bI_1} (A) \cong L_{\bI_2} \Gamma_{\bI_1} (A) \cong L_{\bI_2} (A)$ in the Verdier quotient $\widehat \bK / \Loc(\bI_1 \cap \bI_2)$; in this Verdier quotient, hence have the distinguished triangle
\[
L_{\bI_2} (A) \to A \to L_{\bI_1} (A) \to \Sigma L_{\bI_2} (A).
\]
But since $L_{\bI_1} (A) \cong \Gamma_{\bI_2} (A)$, it is in $\Loc(\bI_2)$, and $\Sigma L_{\bI_2} (A)$ is in $\Loc(\bI_2)^{\perp}$, and so the map $L_{\bI_1} (A) \to \Sigma L_{\bI_2} (A)$ is 0. By the splitting lemma for triangulated categories, one has an isomorphism of distinguished triangles
\begin{center}
\begin{tikzcd}
L_{\bI_2} (A) \arrow[r] \arrow[d, no head, Rightarrow, no head] & A  \arrow[r] \arrow[d, "\cong"]      & L_{\bI_1} (A) \arrow[r, "0"] \arrow[d, no head, Rightarrow, no head] & \Sigma L_{\bI_2} (A) \\
L_{\bI_2} (A)  \arrow[r]                               & L_{\bI_2} (A) \oplus L_{\bI_1} (A) \arrow[r] & L_{\bI_1} (A) \arrow[r]                                     & \Sigma L_{\bI_2} (A).
\end{tikzcd}
\end{center}
The theorem follows.
\end{proof}

\subsection{A generalized Carlson's theorem} We now prove the generalized Carlson theorem for M$\Delta$Cs.

\bth{carlson-gen} Assume that $\bK$ is a small M$\Delta$C which is the compact part of an associated big M$\Delta$C, $\widehat \bK$, and that $\bK$ satisfies either 
\begin{enumerate}
    \item[(i)] $\bK$ is generated by an object $G$ as a thick subcategory, or
    \item[(ii)] $\Spc \bK$ is Noetherian.   
\end{enumerate}
Then for any indecomposable object $A \in \bK$, the closed set $V(A)$ is connected in $\Spc \bK$.
\eth

\begin{proof}
Suppose $V(A)$ is disconnected, via 
\[
V:=Z_1 \sqcup Z_2.
\]
Just as in the proof of \thref{coprime-complem-equiv}, if $\Spc \bK$ is Noetherian, then by \cite[Corollary 5.5]{Rowe}, $Z_1=V(B)$ and $Z_2=V(C)$ for some objects $B$ and $C$ in $\bK$, and if $\bK$ has a generator $G$, then by \coref{va-disconn}, we again see that $Z_1=V(B)$ and $Z_2=V(C)$ for some objects $B$ and $C$. It follows that 
\[
V(B \oplus C) = V(B) \cup V(C) = V(A),  
\]
which by \cite[Theorem 4.9]{Rowe} in the Noetherian case and \prref{bij} in the finite generation case implies $\langle B, C \rangle = \langle A \rangle$. Furthermore, 
\[
\Phi( \langle B \rangle  \cap \langle C \rangle)= V(B) \cap V(C) = \varnothing.
\]
By \leref{decomp-quot}, $A \cong L_{\bI_1} A \oplus L_{\bI_2} A \in \widehat \bK$. Thus $A$ is not indecomposable. 
\end{proof}
\subsection{A Carlson's theorem for the central cohomological support} Recall that 
the cohomology ring of an M$\Delta$C, $\bK$, is the ring $\End^\bullet_{\bK} (\unit)$. 

Given a collection ${\mathcal K}$ of objects that generates the M$\Delta$C, $\bK$, as a triangulated category, the categorical center $C_{\bK}^\bullet$ of the cohomology ring of $\bK$ is defined to be the graded subring spanned by all $g \in \Hom_{\bK}(\unit, \Sigma^i \unit)$, such that for every $M\in {\mathcal K}$, the following diagram commutes, where the isomorphisms are structure isomorphisms for an M$\Delta$C:
\[
\begin{tikzcd}
\unit \otimes M \arrow[d, "g \otimes \id_M"] \arrow[r, "\cong"] & M          & M \otimes \unit \arrow[l, "\cong"] \arrow[d, "\id_M \otimes g"] \\
\Sigma^i \unit \otimes M \arrow[r, "\cong"]                     & \Sigma^i M & M \otimes \Sigma^i \unit \arrow[l, "\cong"]  
\end{tikzcd}
\]
see \cite[Section 1.4]{NVY3}.
The {\it{central cohomological support variety}} of an object $A \in \bK$ is then defined to be
\[
W_C(A):=\{\mf{p} \in \Proj C^\bullet_{\bK} : I(A) \subseteq \mf{p}\},
\]
where $I(A)$ is the annihilator ideal of $\End^\bullet_{\bK} (A)$ in $C^\bullet_{\bK}$ for the action
\[
g \cdot h = \Sigma^i(h) (g \otimes \id_A), 
\quad \mbox{for} \; \; g \in \Hom_{\bK}(\unit, \Sigma^i \unit), h \in \Hom_{\bK}(A, \Sigma^j A),
\]
where we use the natural isomorphism $(\Sigma^i \unit) \otimes A \cong \Sigma^i A$. We say that 
$\bK$ satisfies the weak finite-generation property if for all objects $A \in \bK$, $\End^\bullet_{\bK} (A)$ is a finitely generated $C_{\bK}^\bullet$-module, see
\cite[Section 1.5]{NVY3} for details. When this is satisfied, we have a continuous map
\[
\Spc \bK \xrightarrow{\rho} \Proj C^\bullet_{\bK}
\sqcup \{ \mbox{irrelevant ideal} \}
\]
which is surjective onto $\Proj C^\bullet_{\bK}$,
see \cite[Theorem 7.2.1]{NVY3}. Assume additionally that $W_C$ satisfies the weak tensor product property, that is,
\[
\bigcup_{B \in \bK} W_C(A \otimes B \otimes C) = W_C(A) \cap W_C(C).
\]
Then we have a map
\[
\Proj C^\bullet_{\bK} \xrightarrow{\eta} \Spc \bK
\]
such that $\rho \eta = \id$, see \cite[Proposition 8.1.1]{NVY3}. We obtain a version of the Carlson theorem for central cohomological support.

\bco{Carlson-cohom}
Let $\bK$ be a small M$\Delta$C satisfying either
\begin{enumerate}
    \item[(i)] $\bK$ is generated by an object $G$ as a thick subcategory, or
    \item[(ii)] $\Spc \bK$ is Noetherian.
\end{enumerate} Suppose that $\bK$ also satisfies weak finite generation, and that the central cohomological support satisfies the weak tensor product property. If $A$ is an indecomposable object, then $W_C(A)$ is connected.
\eco

\begin{proof}
We claim that $\rho(V(A))=W_C(A)$: on one hand, we have $\rho(V(A)) \subseteq W_C(A)$ by \cite[Corollary 6.2.5]{NVY3}, and on the other hand, if $\mf{p} \in W_C(A)$, then $\eta(\mf{p}) \in V(A)$ by \cite[Section 8.1.]{NVY3}, and $\rho(\eta(\mf{p}))=\mf{p}$. By \thref{carlson-gen}, $V(A)$ is connected; since $W_C(A)$ is a continuous image of $V(A)$, it follows that $W_C(A)$ is connected as well. 
\end{proof}

Now let $\bC$ be a finite tensor category (c.f.~ \cite{EO}), e.g. $\bC=\modd(H)$ for $H$ a finite-dimensional Hopf algebra. Set $\bK=\ul{\bC}$ the stable category of $\bC$, and set $\widehat \bK = \ul{\Ind(\bC)}$; then $\widehat \bK$ is a compactly-generated M$\Delta$C which has $\widehat \bK^c=\bK$, see \cite[Theorem A.0.1]{NVY3}. The cohomology ring of $\ul{\bC}$ is canonically isomorphic to the cohomology ring of $\bC$. In the definition of categorical center of the cohomology ring of $\ul{\bC}$, one takes the collection ${\mathcal K}$ to be the (finite) set 
of simple objects of $\bC$. Moreover, $\bK$ is generated by a single object $G$ (the direct sum of the simple objects of $\bC$) as a thick subcategory.

From Corollary \ref{cCarlson-cohom}, we immediately obtain the following:

\bco{Carlson-cohom-2}
Let $\bC$ be a finite tensor category such that $\bK:=\ul{\bC}$ satisfies the weak finite generation condition, and such that the central cohomological support satisfies the weak tensor product property. If $A$ is indecomposable, then $W_C(A)$ is connected.
\eco

Bergh, Plavnik and Witherspoon  \cite[Corollary 6.4]{BPW2021} proved a version of Carlson's Connectedness Theorem for cohomological support for finite tensor categories satisfying the Etingof--Ostrik Finite Generation Conjecture \cite{EO}.
Under weak finite generation and the weak tensor product property, Theorem \ref{tcarlson-gen} recovers the analogue of their theorem for the central cohomological support. We note that the validity of the Etingof--Ostrik Finite Generation Conjecture \cite{EO} for the Drinfeld center of $\bC$ (a braided finite tensor category) implies the weak finite generation condition for $\ul{\bC}$, and in fact, the latter condition is a consequence of just one of the two conditions of the former conjecture.

A Carlson Connectedness Theorem in a more general framework that cohomological supports was obtained by Buan, Krause, Snashall and Solberg,\cite[Proposition 6.6]{BKSS} but it relies on a long list of technical conditions.
\section{Examples}
\label{ex}
In this section we illustrate our results with examples from algebraic geometry and group representation theory.
\subsection{Derived categories of schemes}
Let $X$ be a scheme. Denote by $D(X)$ the derived category of coherent sheaves on $X$, and $D^{\perf}(X)$ the perfect derived category of $X$. We recall several properties of $D(X)$ and $D^{\perf} (X).$
\begin{enumerate}
\item[(i)] $D^{\perf}(X)$ is the compact part of $D(X)$ \cite[Lemma 3.5]{Rouquier2010}. 
\item[(ii)] There is a bijection between Thomason closed subsets of $X$ and thick ideals of $D^{\perf}(X)$ \cite[Theorem 3.15]{Thomason1997}. We denote the thick ideal associated to $Z$ by $\bI(Z)$. 
\item[(iii)] $D(X) / \Loc(\bI(Z)) \cong D(X \backslash Z)$ \cite[Lemma 3.4 and Lemma 3.10]{Rouquier2010}. 
\end{enumerate}

Suppose $X$ is a disconnected scheme with $X= X_1 \sqcup X_2$ with $X_1$ and $X_2$ closed. Denote by $\bI_i$ the ideal of objects with (universal) support contained in $X_i$. Set $\bL_i:=\Loc(\bI_i)$. Then \thref{CRT} implies that
\begin{align*}
D(X ) & \cong D(X) / \Loc(\bI_1 \cap \bI_2)\\
&\cong D(X)/\bL_1 \oplus D(X)/ \bL_2  \\
&\cong D(X_2) \oplus D(X_1).
\end{align*}

\subsection{Stable module categories of finite groups}

Let $\bK=\underline{\modd}(\kk G)$ for $G$ a finite group and $\kk$ a field of characteristic dividing the order of $G$, and $\widehat \bK = \underline{\Mod}(\kk G).$ The dimension of the Balmer spectrum is equal to the $p$-rank $r$ of $G$, by Quillen's Stratification Theorem \cite{Quillen} and the fact that the Balmer spectrum of $\bK$ is homeomorphic to the $\Proj$ of the cohomology ring of $G$ (\cite[Theorem 3.4]{BCR1997}, \cite[Corollary 5.10]{Balmer}). Denote the irreducible components of $\Spc \bK$ by $Z_1,..., Z_m$, and denote 
\[
\bI : = \{ A \in \bK : \dim V(A) <r\}
\]
and 
\[
\bJ :=\{A \in \bK: V(A) \textrm{ does not contain any } Z_i \textrm{ for }i=1,...,m\}.
\]

Using the results of Carlson--Donovan--Wheeler and Carlson \cite{CDW1994, Carlson1996}, Benson described the quotient $(\bK/\bI)^{\idem}$ \cite[Theorem 2.1]{Benson1996}. Furthermore, Benson and Krause described the quotient $(\bK/\bJ)^{\idem}$ \cite[Proposition 4.1]{BK2000}. Denote by $\mc{M}$ the collection of all maximal ideals of $\bK$ and $\mc{M}'$ the collection of maximal ideals of $\bK$ so that the corresponding irreducible component of $\Spc \bK$ is of dimension $r$. We observe that
\[
\bI = \bigcap_{\bM \in \mc{M}'} \bM
\]
and
\[
\bJ = \bigcap_{\bM \in \mc{M}} \bM.
\]
It follows that the Decomposition \thref{CRT} generalizes the results of Benson \cite{Benson1996} and Benson--Krause \cite{BK2000}.

\appendix
\section{Correspondence between thick semiprime ideals and Thomason closed sets in Balmer spectra}
\label{top-balmer-ideals}

\subsection{}  We refer to a {\it noncommutative-multiplicative} (or nc-multiplicative) collection of objects as a collection of objects $\mc{M}$ such that if $A$ and $B$ are in $\mc{M}$, there is an object $C$ such that $A \otimes C \otimes B \in \mc{M}$. 

To identify the quasicompact open subsets of the Balmer spectrum, we first need to state a lemma. The lemma below for nc-multiplicative collections was originally stated and proved for multiplicative collections (that is, collections of objects $\mc{M}$ where $A$ and $B$ in $\mc{M}$ imply that $A \otimes B \in \mc{M}$) in \cite[Theorem 3.2.3]{NVY1}. The same proof works in this more general setting. 

\ble{nc-mult}
Let $\mc{M}$ be an nc-multiplicative collection of objects of a monoidal triangulated category $\bK$. Let $\bI$ be a proper thick ideal of $\bK$ that intersects $\mc{M}$ trivially. Then if $\bP$ is a maximal element of the set
\[
X(\mc{M}, \bI) :=\{ \bJ \text{ a thick ideal of }\bK: \bJ \supseteq \bI, \bJ \cap \mc{M}= \varnothing \},
\]
at least one of which exists by Zorn's lemma, then $\bP$ is a prime ideal of $\bK$. 
\ele

\subsection{} Let $\bK$ be an M$\Delta$C. Our goal is to establish a correspondence between semiprime ideals of $\bK$ and Thomason closed sets of the Balmer spectrum of $\bK$, under the assumption that $\bK$ has a generator, that is, there exists an object $G$ which generates $\bK$ as a thick subcategory. In particular, this applies to stable categories of finite tensor categories, where the generator is the direct sum of simple objects. This generalizes the foundational classification theorem of Balmer \cite[Theorem 4.10]{Balmer}. Under the assumption that $\Spc \bK$ is Noetherian, these results were established by Rowe \cite{Rowe}. When all prime ideals of $\bK$ are completely prime, similar results were obtained by Mallick and Ray \cite{MallickRay}. A version of the classification theorem was also stated in \cite[Proposition 5.1]{BKS1} using the language of ideal lattices.

\begin{assumption}
\label{gen-assumption}
Throughout the rest of this appendix, we assume that $\bK$ is an M$\Delta$C with a generator $G$. 
\end{assumption}

\subsection{} The following result shows that we can represent a closed set as the universal support of an object. 

\ble{union-vs-fin}
For any finite collection $\mc{S}$ of objects of $\bK$, there exists an object $A$ such that $V(\mc{S})= V(A)$.
\ele

\begin{proof}
It is clear from the definitions that
\[
V(A) \cap V(B) = \bigcup_{C \in \bK} V(A \otimes C \otimes B)=V(A \otimes G \otimes B).
\]
Since 
\[
V(\mc{S}) = \bigcap_{A \in \mc{S}} V(A),
\]
and since $\mc{S}$ is finite, say $\mc{S}=\{A_1,..., A_n\}$, by induction we obtain
\[
V(\mc{S}) = V(A_1 \otimes G \otimes A_2 \otimes G \otimes... \otimes A_n) .
\]
\end{proof}

A version of the following lemma in the language of ideal lattices can be found in Buan--Krause--Solberg \cite[Lemma 3.2]{BKS1}.

\ble{qc-open}
The quasicompact open sets of $\Spc \bK$ are precisely the complements $V(A)^c$ of the universal support of $A$, for objects $A \in \bK$.
\ele

\begin{proof}
We first show that $V(A)^c$ is quasicompact, for any object $A$. Suppose 
\[
V(A)^c = \bigcup_{i \in I} V(\mc{S}_i)^c
\]
for some collections of objects $\mc{S}_i$ over $i \in I$. Set $\mc{S} = \bigcup_{i \in I} \mc{S}_i$, and note that 
\[
V(A)^c = V(\mc{S})^c.
\]
We claim that we can pick a finite collection $B_1,..., B_n \in \mc{S}$ such that 
\[
B_1 \otimes \bK \otimes B_2 \otimes .... \otimes \bK \otimes B_n \subseteq \langle A \rangle.
\]
Equivalently, since $G$ generates $\bK$ as a thick subcategory, we claim that we can pick $B_1,..., B_n \in \mc{S}$ such that
\[
B_1 \otimes G \otimes B_2 \otimes... \otimes G \otimes B_n \subseteq \langle A \rangle.
\]
Assume to the contrary, that no such sequence is in $\langle A \rangle$. Then we can define $\mc{M}$ as the collection of all objects of the form $B_1 \otimes G \otimes B_2 \otimes... \otimes G \otimes B_n$, over all $B_i \in \mc{S}$. The collection $\mc{M}$ is clearly an nc-multiplicative collection, and intersects $\langle A \rangle$ trivially by assumption. Therefore, by \leref{nc-mult}, we can pick a prime ideal $\bP$ which contains $\langle A \rangle$, and is disjoint from $\mc{M}$ (and in particular is disjoint from $\mc{S}$). But then $\bP$ is in $V(A)^c$, and is not in $V(\mc{S})^c$, which is a contradiction. 

Therefore, we can indeed pick $B_1,..., B_n \in \mc{S}$ with 
\[
B_1 \otimes G \otimes B_2 \otimes... \otimes G \otimes B_n \subseteq \langle A \rangle.
\]
Let $B_i \in \mc{S}_i$. Then we claim
\[
V(A)^c = \bigcup_{i=1}^n V(\mc{S}_i)^c.
\]
The containment $\supseteq$ is clear by assumption. For $\subseteq,$ suppose $\bP$ is in $V(A)^c$, in other words, $A \in \bP$. Then $B_1 \otimes G \otimes B_2 \otimes... \otimes B_n \in \bP$, and thus $B_1 \otimes \bK \otimes B_2 \otimes \bK \otimes... \otimes B_n \subseteq \bP$, which implies by primeness that some $B_i \in \bP$. It follows immediately that $\bP \in V(\mc{S}_i)^c$. This completes the proof that $V(A)^c$ is quasicompact.

For the other direction, suppose we have a quasicompact open set $V(\mc{S})^c$. Note that 
\[
V(\mc{S})^c = \bigcup_{A \in \mc{S}} V(A)^c.
\]
By quasicompactness, there exists a finite subset $\mc{S}' \subseteq \mc{S}$ with 
\[
V(\mc{S})^c = \bigcup_{A \in \mc{S}'} V(A)^c = V(\mc{S}')^c.
\]
The result then follows from \leref{union-vs-fin}.
\end{proof}

\subsection{} The following is now an immediate corollary of \leref{qc-open}, since a topological space is Noetherian if and only if every open subset is quasicompact. In the commutative situation, this was proven in \cite[Corollary 2.17]{Balmer}.

\bco{noeth-va}
$\Spc \bK$ is Noetherian if and only if every closed subset of $\Spc \bK$ is equal to $V(A)$ for some object $A \in \bK$.
\eco

We have an additional useful corollary of \leref{qc-open}.

\bco{va-disconn}
If $A$ is an object of $\bK$ such that $V(A)=V(\mc{S}_1) \sqcup V(\mc{S}_2)$ for some collections of objects $\mc{S}_1$ and $\mc{S}_2$, then there exist objects $A_1$ and $A_2$ in $\bK$ with $V(\mc{S}_i)= V(A_i)$ for $i=1,2$.
\eco

\begin{proof}
    Set $U_i=V(\mc{S}_i)^c$ for $i=1,2$. Then $U_1 \cup U_2 = \Spc \bK$, and $U_1 \cap U_2 = V(A)^c$. Hence both $U_1 \cup U_2$ and $U_1 \cap U_2$ are quasicompact. This implies that both $U_1$ and $U_2$ are quasicompact, by point-set topology. 
\end{proof}

\subsection{} We are now able to give a characterization of the Noetherian property of $\Spc \bK$, in the case that $\bK$ has a finite generating set and either left or right duals. This contrasts the analogous case for commutative rings, where the Noetherianity of the prime spectrum of a ring cannot be used to make any implication about chain conditions for ideals.

\bpr{noeth-rigid}
Suppose $\bK$ is left or right rigid. Then the following are equivalent:
\begin{enumerate}
\item[(a)] $\Spc \bK$ is Noetherian;
\item[(b)] finitely-generated thick ideals of $\bK$ satisfy the descending chain condition.
\end{enumerate}
\epr

\begin{proof}
It is straightforward that if $\Spc \bK$ is Noetherian, then finitely-generated thick ideals satisfy the descending chain condition, see \cite[Lemma B.0.1]{NVY3}. 

For the other direction, suppose $\bK$ satisfies (b). By \coref{noeth-va}, we must just show that a closed subset $V(\mc{S})$ is equal to $V(A)$ for some object $A$, for any collection $\mc{S}$ of objects of $\bK$. Pick an element $A_1 \in \mc{S}$, and set $A(1):=A_1$. We now inductively define an object $A(n)$ in the following way. Suppose there is some object $A_{n} \in \mc{S}$ so that $\langle A(n-1) \otimes G \otimes A_n \rangle$ is properly contained in $\langle A(n-1)\rangle$. Then we define $A(n):=A(n-1) \otimes G \otimes A_n$. In other words, we have inductively
\[
A(n) = A_1 \otimes G \otimes A_2 \otimes G \otimes... \otimes G \otimes A_n,
\]
and a strictly decreasing chain of thick ideals
\[
\langle A(1) \rangle \supsetneq \langle A(2) \rangle \supsetneq \langle A(3) \rangle \supsetneq ...
\]
It is clear that if we can continue this process for all $n \geq 1$, then we obtain an infinite descending chain of principal thick ideals, which contradicts our assumption that (b) holds. Therefore, at some $n$, this process must terminate: we have constructed the object $A(n)$ so that for any additional $A_{n+1} \in \mc{S}$, we have
\[
\langle A(n) \rangle = \langle A(n) \otimes G \otimes A_{n+1} \rangle.
\]

But now we claim
\[
V(\mc{S})=V(A(n)).
\]
Indeed, to show $V(\mc{S})\subseteq V(A(n))$, if $\bP$ intersects $\mc{S}$ trivially, then it clearly also does not contain $A(n)$, since if a prime contains $A(n)$ then it must contain some $A_i$ (which are elements of $\mc{S}$). For $V(A(n))\subseteq V(\mc{S})$, if $\bP$ does not contain $A(n)$, then it does not contain $A(n) \otimes G \otimes A_{n+1}$ for any $A_{n+1} \in \mc{S}$, since $A(n)$ is in the thick ideal generated by $A(n) \otimes G \otimes A_{n+1}$. Hence, $\bP$ does not contain any element of $\mc{S}$, and so $\bP \in V(\mc{S})$. 

Since we have produced a single object $A(n)$ so that $V(\mc{S})$ can be given as $V(A(n))$ for any collection $\mc{S}$ of objects of $\bK$, we conclude that $\Spc \bK$ is Noetherian. 
\end{proof}

\subsection{} Recall that a {\it Thomason subset} of a topological space is a subset $S$ of the form
\[
S=\bigcup_{i \in I} S_i,
\]
where each $S_i$ is closed and has quasicompact complement. We now identify the Thomason subsets of the Balmer spectrum, using \leref{qc-open}.

\bco{thomason-balmer}
The Thomason subsets of $\Spc \bK$ are precisely the subsets of the form
\[
\bigcup_{i \in I} V(A_i),
\]
over all collections of objects $\{A_i\}_{i \in I}$ for any index set $I$. When $\Spc \bK$ is Noetherian, these coincide with the specialization-closed subsets.
\eco

\begin{proof}
This follows directly from \leref{qc-open}.
\end{proof}

\subsection{} The following is a noncommutative version of \cite[Theorem 4.10]{Balmer}. 

\bpr{bij}
The maps
\begin{enumerate}
\item[(a)] $\Phi$ from thick ideals of $\bK$ to specialization-closed subsets of $\Spc \bK$, defined by 
\[
\bI \mapsto \{\bP \in \Spc \bK : \bI \not \subseteq \bP\}=\bigcup_{A \in \bI} V(A),
\]
\item[(b)] and $\Theta$, from specialization-closed subsets of $\Spc \bK$ to thick ideals of $\bK$, defined by
\[
S \mapsto \bigcap_{\bP \not \in S}\bP
\]
\end{enumerate}
descend to an order-preserving bijection between semiprime thick ideals of $\bK$ and Thomason subsets of $\Spc \bK$.
\epr

\begin{proof}


Note that $\Phi(\bI)$ is a Thomason subset of $\Spc \bK$ for any thick ideal $\bI$ by \coref{thomason-balmer}, and that $\Theta(S)$ is a semiprime thick ideal, for any Thomason subset $S$ of $\Spc \bK$. It is clear that both maps are order preserving: if $\bI \subseteq \bJ$ are thick ideals of $\bK$, then 
\[
\{\bP \in \Spc \bK : \bI \not \subseteq \bP\} \subseteq \{\bP \in \Spc \bK : \bJ \not \subseteq \bP\},
\]
and if $S \subseteq T$ are Thomason subsets of $\Spc \bK$, then 
\[
\bigcap_{\bP \not \in S} \bP\subseteq \bigcap_{\bQ \not \in T} \bQ.
\]



We now check that the two maps described are inverse to one another. Let $\bI$ be semiprime. We compute $\Theta( \Phi(\bI))$:
\begin{align*}
\bI \mapsto \{\bP : \bI \not \subseteq \bP\} & \mapsto \bigcap_{\bQ \not \in \{\bP : \bI \not \subseteq \bP\}} \bQ \\
&=\bigcap_{\bI \subseteq \bP} \bP\\
&= \bI.
\end{align*}
The last step uses the semiprimeness assumption. 

On the other hand, to compute $\Phi(\Theta(S))$ for a Thomason subset $S$, we find
\begin{align*}
S  \mapsto \bigcap_{\bP \not \in S} \bP &\mapsto \left \{\bQ : \bigcap_{\bP \not \in S} \bP \not \subseteq \bQ \right \}.
\end{align*}
It is clear by definition that the right hand side of this equation, that is, $\Phi(\Theta(S))$, is contained in $S$. To obtain equality, we must appeal to the assumption that $S$ is Thomason; in other words, using \coref{thomason-balmer}, we can write $S= \bigcup_{i \in I} V(A_i)$ for some collection of objects $A_i$. Suppose $\bQ \in S$, and hence, $\bQ \in V(A_i)$ for some $i \in I$. We now just note that if $\bP \not \in S$, we have $A_i \in \bP$, since otherwise $\bP$ would be in $V(A_i) \subseteq S$. Therefore, if $\bP \not \in S$, we have $\bP \not \subseteq \bQ$, given that $\bP$ contains $A_i$ and $\bQ$ does not. This shows $\bQ$ is in $\Phi(\Theta(S))$, and we are done.
\end{proof}

\bpr{fin-gen-ideals}
Let $\bK$ be left or right rigid M$\Delta$C. Then the bijection given in \prref{bij} descends further to a bijection between finitely-generated ideals of $\bK$ and closed Thomason subsets of $\Spc \bK$. In other words, an ideal $\bI$ is finitely-generated if and only if $\Phi(\bI)$ is closed. 
\epr

\begin{proof}
It is clear that if $\bI$ is finitely generated, say by $A_1,..., A_n$, then it is straightforward that 
\[
\Phi(\bI) = \bigcup_{i=1}^n V(A_i),
\]
a closed set. 

On the other hand, suppose $\Phi(\bI)$ is closed. Since $\bK$ is generated by a single object, by \leref{qc-open}, there exists an object $A$ with $\Phi(\bI)=V(A)$. We know that 
\[
\bI = \bigcap_{A \in \bP} \bP= \langle A \rangle,
\]
using both \prref{bij} and the rigidity assumption, since rigidity implies all ideals are semiprime \cite[Proposition 4.1.1]{NVY2}.
\end{proof}



\begin{thebibliography}{99}
\bibitem{Amit1} 
S.A. Amitsur, {\em A general theory of radicals. I: Radicals in complete lattices},  Amer. J. Math. {\bf 74} (1952), 774--786.

\bibitem{Amit2}
S.A. Amitsur, {\em A general theory of radicals. II: Radicals in rings and bicategories}, Amer. J. Math. {\bf 75} (1954), 100--125.

\bibitem{Amit3} 
S.A. Amitsur, {\em A general theory of radicals. III: Applications}, Amer. J. Math. {\bf 75} (1954), 126--136.


\bibitem{BS} S. Balchin and G. Stevenson, {\em{Big categories, big spectra}}, 
preprint arXiv:2109.11934.  

\bibitem{Balmer}
P. Balmer, {\em{The spectrum of prime ideals in tensor triangulated categories}},
J. Reine Angew. Math. {\bf{588}} (2005), 149--168.

\bibitem{Balmer2} P. Balmer, {\em{Supports and filtrations in algebraic geometry and modular representation theory}},
Amer. J. Math. {\bf{129}} (2007), no.5, 1227--1250.


\bibitem{BalmerFavi0} P. Balmer and G. Favi, 
{\em{Gluing techniques in triangular geometry}},
Q. J. Math. {\bf{58}} (2007), no.4, 415--441.


\bibitem{BalmerFavi}
P. Balmer and G. Favi,
{\em{Generalized tensor idempotents and the telescope conjecture}},
Proc. Lond. Math. Soc. (3) {\bf 102} (2011), no. 6, 1161--1185.


\bibitem{Balmer-Schlichting}
P. Balmer and M. Schlichting,
{\em{Idempotent completion of triangulated categories}},
J. Algebra {\bf 236} (2001), no. 2, 819--834.

\bibitem{Benson1996}
D. J. Benson, 
{\em{Decomposing the complexity quotient category}},
Math. Proc. Camb. Phil. Soc. {\bf 120} (1996), no. 4, 589--595.

\bibitem{BCR1997}
D. J. Benson, J. F. Carlson, and J. Rickard,
{\em{Thick subcategories of the stable module category}},
Fund. Math. {\bf 153} (1997), no.1, 59–80.

\bibitem{BIK2012}
D. J. Benson, S. B. Iyengar, and H. Krause, {\em{Representations of finite groups: local cohomology and support}}, 
Oberwolfach Seminars vol. 43, Birkh\"auser/Springer Basel AG, Basel, 2012.

\bibitem{BIK2013}
D. Benson, S. B. Iyengar, and H. Krause, 
{\em{Module categories for group algebras over commutative rings}},
with an appendix by G. Stevenson,
J. K-Theory {\bf 11} (2013), no. 2, 297--329.

\bibitem{BK2000}
D. J. Benson and H. Krause, 
{\em{Generic idempotent modules for a finite group}},
special issue dedicated to Klaus Roggenkamp on the occasion of his 60th birthday,
Algebr. Represent. Theory {\bf 3} (2000), no. 4, 337--346.

\bibitem{Benson-Witherspoon} D. J. Benson and S. Witherspoon, 
{\em{Examples of support varieties for Hopf algebras with noncommutative tensor products}}, Arch. Math. (Basel) {\bf{102}} (2014), no. 6, 512--520. 

\bibitem{BPW2021}
P. A. Bergh, J. Y. Plavnik, and S. Witherspoon,
{\em{Support varieties for finite tensor categories: complexity, realization, and connectedness}},
J. Pure Appl. Algebra {\bf 225} (2021), no.9, Paper No. 106705, 21 pp.

\bibitem{BKSS} A.B. Buan, H. Krause, N. Snashall, and {\O}. Solberg, {\em{Support varieties--an axiomatic approach}},  
Math. Z. {\bf{295}} (2020), 395--426. 

\bibitem{BKS1} A. B. Buan, H. Krause, and Ø. Solberg,
{\em{Support varieties: an ideal approach}},
Homology Homotopy Appl. {\bf 9} (2007), no. 1, 45--74.

\bibitem{Carlson1984}
J. F. Carlson,
{\em{The variety of an indecomposable module is connected}},
Invent. Math. {\bf 77} (1984), no. 2, 291--299.

\bibitem{Carlson1996}
J. F. Carlson, 
{\em{The decomposition of the trivial module in the complexity quotient category}},
J. Pure Appl. Algebra {\bf 106} (1996), no. 1, 23--44.

\bibitem{CDW1994}
J. F. Carlson, P. W. Donovan, and W. W. Wheeler, 
{\em{Complexity and quotient categories for group algebras}}, 
J. Pure Appl. Algebra {\bf 93} (1994), no. 2, 147--167.

\bibitem{EGNO} P. Etingof, S. Gelaki, D. Nikshych, and V. Ostrik,
{\em{Tensor categories}}, Math. Surveys and Monographs vol. 205,
Amer. Math. Soc., Providence, RI, 2015.

\bibitem{EO} P. Etingof and V. Ostrik,
{\em{Finite tensor categories}},
Mosc. Math. J. {\bf{4}} (2004), no. 3, 627--654.

\bibitem{GW} K. R. Goodearl and R. B.  Warfield, Jr.
{\em{An introduction to noncommutative Noetherian rings}}, 
Second ed., London Math. Soc. Stud. Texts {\bf{61}},
Cambridge Univ. Press, Cambridge, 2004.

\bibitem{Jac} N. Jacobson, {\em Basic Algebra I}, W.H. Freeman and Company, 2nd edition, 1985. 

\bibitem{Krause0} H. Krause,
{\em{Decomposing thick subcategories of the stable module category}},
Math. Ann. {\bf{313}} (1999), no.1, 95--108.


\bibitem{Krause}
H. Krause,
{\em{Localization theory for triangulated categories}}, in: Triangulated categories, pp. 161--235,
London Math. Soc. Lecture Note Ser., 375, Cambridge Univ. Press, Cambridge, 2010.

\bibitem{Krause2}
H. Krause,
{\em{Central support for triangulated categories}},
Int. Math. Res. Notices IMRN {\bf{2023}}, no. 22, 19773--19800.


\bibitem{Kunz} E. Kunz, 
{\em Introduction to commutative algebra and algebraic geometry}, Birkh\"{a}user, Boston, 1985. 

\bibitem{MallickRay}
V. M. Mallick and S. Ray,
{\em{Noncommutative tensor triangulated categories and coherent frames}}, preprint arXiv:2204.08794, to appear in 
C. R. Math. Acad. Sci. Paris.


\bibitem{NVY1}
D. K. Nakano, K. B. Vashaw, and M. T. Yakimov, {\em{Noncommutative tensor triangular geometry}},
Amer. J. Math. {\bf{144}} (2022), no. 6, 1681--1724.

\bibitem{NVY2}
D. K. Nakano, K. B. Vashaw, and M. T. Yakimov,
{\em{Noncommutative tensor triangular geometry and the tensor product property for support maps}},
Int. Math. Res. Notices IMRN {\bf{2022}}, no. 22, 17766--17796.

\bibitem{NVY3}
D. K. Nakano, K. B. Vashaw, and M. T. Yakimov, 
{\it On the spectrum and support theory of a finite tensor category}, 
Math. Ann. doi:10.1007/s00208-023-02759-8.

\bibitem{Neeman1} A. Neeman, {\em{Triangulated categories}}, 
Ann. Math. Studies {\bf{148}}, Princeton Univ. Press, Princeton, NJ, 2001. 

\bibitem{Plavnik-Witherspoon} J. Plavnik and S. Witherspoon, {\em Tensor products and support varieties for some noncocommutative Hopf algebras}, Algebr. Represent. Theory {\bf{21}} (2018), 259--276. 

\bibitem{Quillen}
D. Quillen,
{\em{The spectrum of an equivariant cohomology ring. I, II}},
Ann. of Math. (2) {\bf 94} (1971), 549–572; ibid. (2) {\bf 94} (1971), 573--602.

\bibitem{Rickard1}
J. Rickard, 
{\em{Idempotent modules in the stable category}},
J. London Math. Soc. (2) {\bf 56} (1997), no. 1, 149–170.


\bibitem{Rowe}
J. Rowe,
{\em{Noncommutative tensor triangulated geometry:
classification via noetherian spectra}}, preprint arXiv:2308.14661.

\bibitem{Rouquier2010}
R. Rouquier,
{\em{Derived categories and algebraic geometry}}, in: Triangulated categories, pp. 351–-370,
London Math. Soc. Lecture Note Ser., 375, Cambridge Univ. Press, Cambridge, 2010.

\bibitem{Stevenson}
G. Stevenson,
{\em{A tour of support theory for triangulated categories through tensor triangular geometry}},
Building bridges between algebra and topology, 63–101,
Adv. Courses Math. CRM Barcelona, Birkhäuser/Springer, Cham, 2018.

\bibitem{SA} M. Suarez-Alvarez, {\em{The Hilton--Heckmann argument for the anti-commutativity of cup products}}, Proc. Amer. Math. Soc. {\bf{132}} (2004), no. 8, 2241--2246.

\bibitem{Thomason1997}
R. W. Thomason,
{\em{The classification of triangulated subcategories}},
Compositio Math. {\bf{105}} (1997), no. 1, 1--27.
\end{thebibliography}
\end{document}